\documentclass[12pts]{amsart}
\usepackage{amssymb}
\usepackage{epsfig}

\newtheorem{thm}{Theorem}[section]
\newtheorem{lem}[thm]{Lemma}

\newtheorem{prop}[thm]{Proposition}

\theoremstyle{remark}

\newtheorem{exmp}[thm]{Example}

\theoremstyle{definition}
\newtheorem{defn}[thm]{Definition}

\newcommand{\bZ}{\mathbb{Z}}
\newcommand{\bC}{\mathbb{C}}
\newcommand{\bR}{\mathbb{R}}
\newcommand{\hs}{{\mathbb{H}}^3}
\newcommand{\bU}{\mathcal{U}}
\newcommand{\cV}{\mathcal{V}}
\newcommand{\cL}{\mathcal{L}}
\newcommand{\cM}{\mathcal{M}}
\newcommand{\cC}{\mathcal{C}}
\newcommand{\pslc}{\mathrm{PSL}_2(\mathbb{C})}
\newcommand{\slc}{\mathrm{SL}_2(\mathbb{C})}
\newcommand{\tr}{\mathrm{tr}}
\newcommand{\re}{\mathrm{Re}}
\newcommand{\im}{\mathrm{Im}}
\newcommand{\hatmu}{\hat{\mu}}
\newcommand{\dev}{{\mathrm{dev}}}
\newcommand{\ra}{\rightarrow}
\newcommand{\tri}{\triangle}

\begin{document}

\title[Deformations of hyperbolic cone-manifolds]
{The maximal tube under the deformations of 
a class of $3$-dimensional hyperbolic cone-manifolds}

\author{Suhyoung Choi}
\address{Department of Mathematics\\ 
Seoul National University \\
151--742 Seoul, Korea} 
\email{shchoi@math.snu.ac.kr} 
\author{Jungkeun Lee} 
\address{}
\email{jklee@math.snu.ac.kr} 
\subjclass{Primary 57M50} 
\date{May 10, 2004} 
\keywords{hyperbolic manifold, cone-manifold, deformations} 
\thanks{The first author gratefully acknowledges  
support from Korea Research Foundation Grant (KRF-2002-070-C00010).}

\begin{abstract}
Recently, Hodgson and Kerckhoff found a small bound on  
Dehn surgered $3$-manifolds from hyperbolic knots 
not admitting hyperbolic structures using deformations of 
hyperbolic cone-manifolds. They asked whether 
the area normalized meridian length squared of 
maximal tubular neighborhoods of 
the singular locus of the cone-manifold
is decreasing and that summed with the cone angle squared
is increasing as we deform the cone-angles.    
We confirm this near $0$
cone-angles for an infinite family of hyperbolic cone-manifolds 
obtained by Dehn surgeries along the Whitehead link complements. 
The basic method is based on explicit holonomy computations using 
the $A$-polynomials and finding the maximal tubes.
One of the key tool is the Taylor expression of a geometric
component of the zero set of the $A$-polynomial 
in terms of the cone-angles. 
We also show a sequence of 
Taylor expressions for Dehn surgered manifolds 
converges to one for the limit hyperbolic manifold. 
\end{abstract} 

\maketitle

\section{Introduction}\label{sec:intro}

%
Recently it was shown by Hodgson and
Kerckhoff(\cite{HK2}) that there is a small universal bound for the
number of nonhyperbolic Dehn-fillings on a hyperbolic manifold
with single cusp. Their argument involves analysis of the
variation of maximal tubes around singularities in cone manifolds
of fixed topological type when the cone angle increases starting
from 0.

Let $M$ be an orientable
3-manifold which admits a complete hyperbolic structure of finite
volume with single cusp. For each slope $\gamma$ of the cusp of
$M$, let $M(\gamma)$ be the manifold obtained by Dehn-filling $M$
along $\gamma$. By a theorem of Gromov and Thurston (\cite{BH}),
$M(\gamma)$ admits a negatively curved metric if the length of the
shortest curve on the boundary of a horoball neighborhood of the
cusp isotopic to $\gamma$ is greater than $2\pi$. But it is not
known whether $M(\gamma)$ admits a hyperbolic structure with the
same hypothesis. 

A hyperbolic cone-manifold of $3$-dimension is a manifold locally modeled 
on open subset of a hyperbolic space or the open region in an open set 
bounded by two totally geodesic planes meeting at a geodesic and two planes
are identified by an elliptic isometry. 
By Thurston's hyperbolic Dehn surgery theorem, if
$\theta > 0$ is small, $M(\gamma)$ for any $\gamma$ admits a hyperbolic cone-structure 
whose singular locus is the added closed curve with some small cone angle $\theta$. 
We denote the resulting cone manifold by $M(\gamma; \theta)$. 
The homotopy class of the singular locus obviously corresponds 
to the closed curve meeting with $\gamma$ once.
(Often $M$ itself will be considered $M(\gamma; 0)$ and 
as having $\theta = 0$.) 


Each of the cone manifold has the maximal tube around its singular locus. If we can
bound from below the radii of the maximal tubes until the cone
angle reaches $2\pi$, then we obtain a nonsingular hyperbolic
structure on $M(\gamma)$. Let $T$ be the flat torus boundary of
the maximal horoball neighborhood of the cusp in $M$. In
\cite{HK2}, Hodgson and Kerckhoff show that if the normalized
length of the geodesic curve on $T$ isotopic to $\gamma$, i.e. the
length measured when the metric  on $T$ is rescaled so that the
area of $T$ is 1, is greater than $7.515$, then we can bound from
below the radii of the maximal tubes in $M_\theta (\gamma)$ until
the cone angle $\theta$ reaches $2 \pi$ and hence $M(\gamma)$ is
hyperbolic. Using this fact, they obtained the universal bound 60
for the number of nonhyperbolic Dehn-fillings on a single-cusped
hyperbolic manifold.  

In one of their lectures, they posed the following question:\\ \\
{\bf Question} (Hodgson and Kerckhoff) 
Let $\{ M_\theta (\gamma): 0 < \theta < \theta_0 \}$ be 
a continuous family of hyperbolic
cone structures on $M(\gamma)$ with singular locus as described
above. Let $\mu = \mu(\theta)$ be the length of the surgery curve
on the boundary $T_\theta $ of the maximal tube around the
singular locus of $M_\theta (\gamma)$ and let $\hatmu$ be the
normalized length $ \mu / \sqrt{{\mathrm {Area}}(T_\theta )}$.
Then are $\hatmu^2 $ and $\hatmu ^2 + \theta ^2 $ decreasing and
increasing
functions of $\theta $ on $[0 , \theta _0 )$, respectively? \\ \\
They showed that if this question had positive answer then we
could also control the radii of the maximal tubes effectively
(see \cite{Ker} for more details).

Our initial result is on the relationship between 
$A$-polynomial and the cone-angle: 

Let $M$ be a hyperbolic manifold of finite volume and two cusps 
and $\{\cM_1 , \cL_1 \}$ be a basis for a peripheral group
group $P$ of the fundamental group of $\pi_1(M)$ corresponding to 
the first cusp of $M$. We choose a basis $\{ \cM_2, \cL_2\}$ 
for the second peripheral group and fix the two bases. 
Let $(p_1, q_1)$ and $(p_2, q_2)$ denote coprime pair of integers.
Let $M(p_1, q_1)$ denote the $3$-manifold obtained from $M$ by 
the $(p_1, q_1)$-Dehn filling the first cusp. 
Let $M(\infty, \infty)(p_2, q_2; \theta)$ denote the hyperbolic cone-manifold 
with a cusp corresponding to the first cusp of $M$ and the second 
cusp has been $(p_2, q_2)$-Dehn filled where the corresponding solid torus
has a cone-type singularity with the cone-angle equal to $\theta$. 
Let $M(p_1, q_1)(p_2, q_2; \theta)$ denote the hyperbolic cone-manifold 
with the first cusp $(p, q)$-Dehn-filled and 
the second cusp has been $(p_2, q_2)$-Dehn-filled with the cone-type
singularity with the cone-angle equal to $\theta > 0$. 

The character variety of $\pslc$-representations of the fundamental group 
of $M$ is the space of characters $\pi_1(M) \ra \bC$ defined by 
taking the traces of holonomies of the fundamental group $\pi_1(M)$ 
(see Culler-Shalen \cite{CS} for details). 

The character variety of $\pslc$-representations of the fundamental group
of $M$ which keep the holonomy of $\cM_1, \cL_1$ parabolic
gives a relation between the eigenvalues $l_1$ and $m_1$ of holonomies of 
$\cL_1$ and $\cM_1$, which is how $A$-polynomials are defined in this paper. 
The so-called geometric component of the zero-locus of the $A$-polynomial is 
a component realized by a deformation of hyperbolic manifold corresponding 
to the cone-structures for the second cusp.  
Using the geometric component of 
the $A$-polynomial together with the Dehn-filling relation
\[p\log{\frac{m}{m_0}}+q\log{\frac{l}{l_0}}=\frac{ \sqrt{-1}\theta }{2},\]
we can obtain Taylor expansions of $m$ and $l$ in terms of
$\theta$, which we need only up to order three. 

We define the $A$-polynomial of $M(p_1, q_1)$ in the same manner 
with respect to $\{\cM_2, \cL_2\}$. 
By the Dehn filling relation, we obtain the Taylor series of $m_2$ and $l_2$ 
corresponding to the geometric component of the $A$-polynomial
as a function of $\theta$. 

Finally, we prove Theorem \ref{thm:convcoeff} showing 
the convergence of the Taylor expansions up to order three 
of $l$ and $m$ in terms of $\theta$ for the $A$-polynomial of $M(p,q)$ to 
the Taylor series of $M(\infty, \infty)$. 
(This will be proved at last Section \ref{sec:convAproof} 
because of the length.)

Our main result is that for an infinite number of of hyperbolic
manifolds $\{W(p_1 , q_1 )\}$ which are obtained from the Whitehead
link complement $W$ by Dehn-fillings on the first torus end, we have a
partial answer to the question of Hodgson-Kerckhoff. 

Let $W$ be the Whitehead link complement, 
and let $\cM _1 , \cL _1 , \cM _2 , \cL _2 $ be
suitably chosen meridians and longitudes for two cusp ends of $W$.

\begin{thm}\label{thm:main} 
Let $\mu = \mu_{p_1, q_1, p_2, q_2} (\theta)$ be
the length of the surgery curve on the boundary of the maximal tube of 
$W(p_1, q_1)(p_2 , q_2; \theta)$ around the singularity. 
Let $\hat{\mu}= \hat{\mu}_{p_1, q_1, p_2, q_2} (\theta)$ 
be the normalized length of the surgery curve. 
If $|p_1| + |q_1|$ is sufficiently large, then for any coprime pair $p_2,q_2$ of
integers except for at most one pair, 
$\hat{\mu}$ is decreasing and $\hat{\mu}^2 + \theta ^2 $
is increasing at $\theta=0$.
\end{thm}

We outline the proof of Theorem. 


For a general hyperbolic manifold $M$ of finite
volume with a distinguished cusp, we detect the maximal horoball
neighborhood of the cusp by finding elements of $\pi_1(M)$ whose
holonomy have the large isometric spheres. 

We use these elements, so called tie classes, 
to find the maximal tube neighborhood of the singularity in
$M(p , q; \theta)$ when $|p| +|q|$ is large and $\theta$ is small.
This follows since the tie classes are stable near $\theta = 0$. 
(See Section \ref{sec:maxtubes} for more details.) 

We now express
the length $\mu_{p,q}(\theta)$ and normalized length $\hatmu_{p,q} (\theta)$
of the surgery curve on the maximal tube around the singularity in
$M_\theta(p , q)$ in terms of the traces of
holonomy of the commutator of a tie class and another element
(see Proposition \ref{prop:formularadius}).

Now we restrict our attention to $W$:
We can exactly compute the holonomy representation of $\pi_1(W)$ corresponding to 
the complete structure.
By looking at the Ford domain, 
$W$ with complete hyperbolic structure decomposes into four ideal tetrahedra.
We can also determine the tie class. 
Each tetrahedron is assigned a complex invariant up to isometry.
These invariants $z_1, z_2, z_3, z_4$ 
satisfy two relations according to two ideal edges of $W$. 
The relation determine a small complex surface parameterizing all hyperbolic 
structures on $W$ near the complete hyperbolic structure. 

Let $m_1, l_1, m_2, l_2$ denote the eigenvalues of the holonomy of 
$\cM_1, \cL_1, \cM_2, \cL_2$ respectively. 
We can write these as functions of the tetrahedral invariants 
$z_1, z_2, z_3, z_4$. 

Next, we compute the holonomy representations as functions of some easily identifiable
variables $x, y$. We easily see that $x = m_1$ and 
we can write $l_1$ as a function of $x$ and $y$. 
Thus, the holonomy representations are functions of $m_1, l_1$. 
From this, we can also find two relations between 
$m_1, l_1$ and $m_2, l_2$.

Using (*), we express $\hatmu_{p_1, q_1, p_2, q_2}^2 (\theta)$ as
\[ k_0 (p_1 , q_1 , p_2 , q_2 ) + k_1 (p_1 , q_1 , p_2 , q_2 )
\theta ^2 + O(\theta ^3),\]
where $k_0$ and $k_1$ are functions
defined for integers $(p_1 , q_1 , p_2 , q_2 )$
with sufficiently large $|p_1|+|q_1|$ and $|p_2|+|q_2|$.
By Theorem \ref{thm:convcoeff}, 
\[k_1(p_1 , q_1 , p_2 , q_2 ) \rightarrow k_1^\infty( p_2 , q_2 )\] 
as $|p_1| + |q_1| \rightarrow \infty$ for some function $k_1^\infty$,
which takes values in some interval $[K_1 , K_2] \subset (-1, 0)$. 

Now, $k_1^\infty(p_2, q_2)$ is the corresponding function for 
$W(\infty, \infty, p_2, q_2)$. We can compute this function 
as in the above general discussion and the main theorem follows.

Dowty also obtained a similar result for figure-eight knot complements
in his doctoral thesis \cite{Dowty} under the supervision of Hodgson. 
Our result generalize his result but we are able to understand 
the effect of Dehn surgery better. 
Our long term hope is that we can answer the question of Hodgson and Kerckhoff for
more general manifolds and with no angle restrictions although there seems to 
be no possible general theory insight.  We think that 
our technique is interesting in that there may be many avenues and 
examples we can consider further and 
serve as a motivation for developing a general theory.

We thank Darryl Cooper, Craig Hodgson, Steven Kerckhoff for many discussions 
and their help. We also thank the Department of Mathematics of 
Stanford University for their great hospitality 
where some of this research was carried out.


\section{Hyperbolic cone-manifolds and hyperbolic Dehn surgery}



In this section, we recall some facts on hyperbolic Dehn surgery theory and 
hyperbolic cone-manifolds. We conclude with a needed result of Neumann and 
Zagier \cite{NZ}.

A hyperbolic manifold is a manifold equipped with a Riemannian
metric whose sectional curvature is the constant $-1$. (In this
paper, we will consider only 3-dimensional
and oriented manifolds.) 
We will denote a simply-connected complete hyperbolic
manifold by $\hs$. The group of orientation preserving isometries
of $\hs$ form a Lie group $\pslc$. If $M$ is a hyperbolic
manifold, each point of $M$ has an open neighborhood isometric to
an open set in $\hs$. We have an isometry $\dev$ from the
universal cover $\tilde M$ onto $\hs$ and a group
homomorphism $\rho : \pi_1 (M) \rightarrow \pslc $ such that 
$\dev \circ \gamma = \rho(\gamma) \circ \dev$ for each 
$\gamma$ in the deck transformation group $\pi_1(M)$. 
Thus, a hyperbolic manifold
is isometric to the quotient space of $\hs$ by the action of a
discrete group of orientation-preserving isometries of $\hs$. 

Hyperbolic cone-manifolds arise in the context of Thurston's hyperbolic Dehn
surgery and are subjects of great interest. (\cite{Ker},
\cite{Ko1}, \cite{Ko2})

\begin{figure}

\begin{center}
\epsfig{file=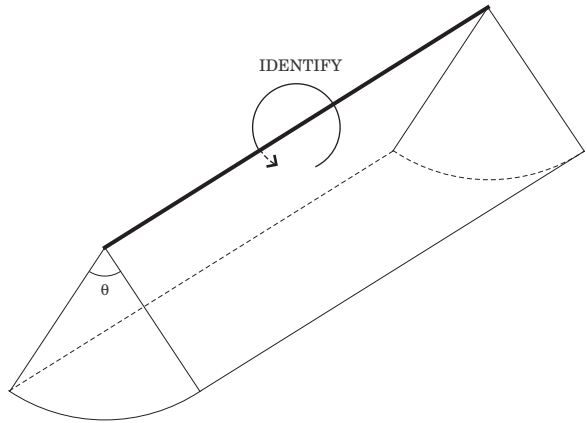, height = 8cm, angle=-90}
\end{center}
\caption{Neighborhood of a singular point in a hyperbolic cone
manifold} \label{fig:conesingularity}
\end{figure}

A 3-manifold $N$ equipped with a metric is a hyperbolic cone-manifold if 
each point of $N$ has an open neighborhood isometric
to an open set in $\hs $ or the quotient metric space 
obtained from an open $3$-ball in $\hs$ 
by removing the domain with boundary in
two totally geodesic planes meeting at a geodesic 
and identifying the two corresponding faces by an isometry fixing the geodesic
(see Figure \ref{fig:conesingularity}).

The set of all points of a hyperbolic cone-manifold $N$ with 
no open neighborhood isometric to an open set in $\hs$ is
called the {\em singular set} or the {\em singularity} of $N$ and is denoted
by $\Sigma _N$( or $\Sigma$ if $N$ is clear from the context). The
singular set is a 1-dimensional submanifold and is a link if $M$
is a closed manifold. To each component of the singular set is
associated the cone angle around the component. 

Let $M$ be a hyperbolic manifold of finite volume with $h$ number of cusps. 
Thurston showed that we can deform $M$ so that $M$ has  
incomplete hyperbolic structures and the metric completions of
some of the deformed structures induces complete hyperbolic
cone-structures on the manifolds obtained by Dehn-fillings on ends of $M$. 
The hyperbolic Dehn surgery theorem states that when we perform
Dehn-fillings on ends of $M$, the resulting manifold admits
complete hyperbolic structures in most cases, i.e., without cone-singularities.


Let $M$ be a hyperbolic manifold of finite volume with $h$
cusps which can be obtained from an ideally triangulated
hyperbolic manifold by Dehn-filling some of the ends. 
Let $\nu$ be
the number of the ideal tetrahedra and let $\cM_1 , \cL_1 ,
\cdots , \cM_h , \cL_h$ be fixed meridians and longitudes for the
ends of $M$. Let the complete hyperbolic structure of $M$
correspond to the point ${\mathbf z}^0 = (z_1^0 , \cdots , z_\nu^0) \in \bC^\nu $. 
Thurston showed that 
the set $\cV$ of points in $\bC^\nu$ near
${\mathbf z}^0$ satisfying certain gluing consistency relations is
a smooth analytic subset of complex dimension $h$ in $\bC^\nu$. 
It was shown that maps ${\mathbf m} = (m_1 , \cdots, m_h): \cV \rightarrow \bC^h$ 
and ${\mathbf l}=(l_1 , \cdots , l_h) : \cV \rightarrow \bC^h$
which assign certain eigenvalues of holonomy
images of $\cM_1 , \cdots ,\cM_h$, $\cL_1 , \cdots ,\cL_h$,
respectively are biholomorphic map at ${\mathbf z}^0$ 
(see Neumann-Zagier \cite{NZ}).

Let $(m_1^0 , \cdots , m_h^0) = {\mathbf m} ({\mathbf z}^0)$ and 
$(l_1^0 , \cdots , l_h^0) = {\mathbf l} ({\mathbf z}^0)$.

\begin{thm}[Neumann-Zagier \cite{NZ}] \label{thm:NZ}
For each $i \in \{1 , \cdots , h \}$, there is a holomorphic
function $\tau_i$ defined on a neighborhood of the origin in
$\bC^h$ such that
\[\log (l_i / l_i ^0) =  
\log (m_i / m_i ^0) \; \tau_i (\log (m_1 / m_1^0), \cdots , \log (m_h / m_h ^0) ). \] 
Moreover for each $i$,
$\tau _i$ is an even function in each variable and 
$\tau _i (0 , \cdots , 0 )$ is the modulus of the flat torus boundary of a cusp
neighborhood for the $i$-th end with respect to $\cM _i , \cL _i$.
In particular each $\tau _i (0 , \cdots , 0 )$ is not a real
number.
\end{thm}
This result will be needed later.

\section{$A$-polynomials and generalized Dehn-fillings}
\label{sec:A-Dehnfilling}

We define the so-called $A$-polynomial for 
multicusped hyperbolic manifolds with respect to a distinguished cusp.
(We may name this relative $A$-polynomial also.) 
We define geometric components of the algebraic set of eigenvalues of 
holonomies, i.e., the zero set of the $A$-polynomial. 
We write the Taylor series of the polynomial relations of geometric components. 
We give some examples. 
Next, we show how to parameterize the geometric component as 
a function of the cone-angles (complex).
Finally, we show that the Taylor series of a geometric component of 
a manifold with one cusp $(p, q)$-Dehn-filled converge to that of 
a manifold without filling as $(p, q) \ra \infty$.

\subsection{Geometric components of the spaces of representation eigenvalues}
\label{subsec:A-polynomial}

Culler-Shalen \cite{CS} defined a character variety of 
a $3$-manifold to be the algebraic set of traces of 
the holonomy of the fundamental group elements ordered in some way. 

We modify Cooper, Culler, Gillet, Long, Shalen \cite{CCGLS}: 
Let $M$ be a hyperbolic $3$-manifold with at least one cusp. 
Fix a cusp and denote by $P$ its fundamental group. 
Denote by $R(\pi_1(M))$ the space of representations of 
the fundamental group $\pi_1(M)$ of $M$ in $\slc$.
We denote by $R(\pi_1(M))_P$ the subset of those representations
whose restrictions to the closed loops in all cusps other than $P$ 
have parabolic or identity images. 
The variety of characters of $\slc$-representations 
of $\pi_1(M)$ is denoted by $X(M)$ and 
$t: R(\pi_1(M)) \ra X(M)$ be the canonical surjective 
projection (see \cite{CS}). We may write $R(M)$ for $R(\pi_1(M))$. 
We denote by $X(M)_P$ the image of $R(M)_P$. 
Using the same reasoning as in \cite{CS}, $X(M)_P$ is 
a dense subset of a finite union of varieties, 
to be denoted by $X''(M)_P$. (We do not claim that they are equal.)

Let $B = \{\cL ,\cM \}$ be a fixed basis for $P$. 
We define the restriction map $r: X''(M)_P \ra X(P)$. 
Define $\tri$ to be the subspace of diagonal representations
in $R(P)$. There is an isomorphism $p_B: \tri \ra \bC^* \times \bC^*$
defined by setting $p_B(\rho) = (l, m)$ if $\rho$ is given by 
\[ \rho(\cL) = \begin{bmatrix} l & 0 \\ 0 & l^{-1} \end{bmatrix} 
\hbox{ and } \rho(\cM) = \begin{bmatrix} m & 0 \\ 0 & m^{-1} \end{bmatrix}. \] 
$t$ induces to a $2-1$-map $t_\tri: \tri \ra X(P)$. 

Denoting $X'(M)_P$ to be the union of irreducible components of $X''(M)_P$ whose 
images under $r$ is complex $1$-dimensional.
For each component $Z'$ of $X'(M)_P$, let $Z$ be the curve 
$t^{-1}_\tri(\overline{r(Z')}) \subset \tri$. 
Define $D_{M, P}$ to be the union of curves $Z$ as $Z'$ varies over 
the components of $X''(M)_P$. 

We say that $D_{M, P}$ is the {\em $A$-set} of $M$ with respect to $P$. 
We note that $A$-set is invariant under the involution
$(l, m) \mapsto (l^{-1}, m^{-1})$.

We define the {\em $A$-polynomial} $A_{M, P}$ of $M$ with respect to $P$ as 
the defining polynomial of the closure of $D_{M, P}$ in $\bC \times \bC$. 

When $M$ has only one cusp, our definition coincide with the definition 
in \cite{CCGLS}.
When obvious, we will drop $P$ from $A_{M, P}$.


\begin{defn}
Let $M$ be a cusped hyperbolic manifold and $B = \{\cL ,\cM \}$ 
be a basis for the fundamental group of a cusp
neighborhood of $M$. 
Suppose that $(l^0 , m^0) \in \bC \times \bC$
equals one of  $(\pm 1,\pm 1)$  and is in the zero set of $A_M (l,m)$. Let
$l(m)$ be a holomorphic function defined on a neighborhood (say $U$) 
of $m^0$ and taking values near $l^0$. We say that the
holomorphic function $l = l(m)$ defined near $(l^0 , m^0) \in \bC \times \bC$ 
is a {\it geometric curve} of the $A$-set at
$(l^0 , m^0)$ if there is a holomorphic family $\{\rho _m  : m \in U \}$ 
of representations of $\pi_1(M)$ into $\slc$ such that each
$\rho_m$ is a lift of a holonomy representation of a hyperbolic
structure on $M$ and
$$\rho _m (\cL) =
\begin{bmatrix} l(m) & * \\ 0 & 1/l(m) \end{bmatrix} , \quad
\rho_m (\cM) =
\begin{bmatrix} m & ** \\ 0 & 1/m \end{bmatrix}.$$
\end{defn}
Clearly, if $l = l(m)$ is a geometric curve of the A-polynomial
$A_M (l , m)$ at $(m^0 , l^0)$, then $A(l(m), m)=0$
for all $m$ near $m^0$. 

\begin{prop}\label{prop:A-set}
A component of $A$-set contains the image of the geometric curve
as a dense set. 
\end{prop}
\begin{proof}
Straightforward.
\end{proof}

A {\em geometric component} is a component of the $A$-set containing 
the geometric curve as a dense set. 
A {\em geometric factor} is a generator of the ideal defining 
the component above.

\begin{exmp} 
Let $M$ be the figure eight knot complement.
$\pi_1(M)$ has a Wirtinger presentation
$$ <\alpha,\beta  : \; \alpha^{-1}\beta \alpha \beta^{-1}
\alpha \beta \alpha^{-1} \beta^{-1} \alpha \beta^{-1}> $$ such
that $\{\alpha , \; \beta ^{-1} \alpha \beta \alpha ^{-2} \beta
\alpha \beta ^{-1} \}$ is a basis for a peripheral subgroup of
$\pi_1 (M)$. The A-polynomial of $M$ with respect to this basis is
$$A_M(l,m) = l m^8 - l m^6 - (l^2 + 2 l + 1) m^4 - l m^2 + l. $$
Note that $A_M(-1 , -1 ) = 0$ and there are two geometric factors
of $A_M(l , m)$ at $(-1 , -1)$ which are
$$l + 1 - 2 \sqrt{-3}(m+1)- (6 + \sqrt{-3})(m+1)^2 - 
( 6 - 2 \sqrt{-3} / 3 )(m+1)^3 + O((m+1)^4)$$
 and 
$$l + 1 + 2 \sqrt{-3}(m+1) - (6 - \sqrt{-3})(m+1)^2 - 
( 6 + 2 \sqrt{-3} / 3 )(m+1)^3 + O((m+1)^4).$$
\end{exmp}

\subsection{Taylor series of geometric curves}
\label{subsec:asymptotic}

Let $M$ be a 3-manifold admitting a complete hyperbolic structure
of finite volume with single cusp. 
Let $\cM , \cL$ be a fixed meridian-longitude pair on
the end of $M$. Let $A_M(l,m)$ be the A-polynomial of $M$ 
with respect to the meridian-longitude pair.

Suppose that deformations of hyperbolic structures
on $M$ near the complete structure gives us 
the following relation for eigenvalues $m$ and $l$ of
$\cM$ and $\cL$ respectively for lifts
of holonomy
representations of nearby hyperbolic structures.
\begin{equation}\label{eq:rel1}
\begin{split}
l = l^0 + a_1 (m - m^0) + 
\frac{a_2}{2} (m -m^0)^2+ \frac{a_3}{6} (m -m^0 )^3 \\+
\mbox{higher order terms}.
\end{split}
\end{equation}
(Here $m^0$ and $l^0$ are the eigenvalues for the lift of a
holonomy representation of the complete structure. Thus each of $m^0$ and
$l^0$ is $\pm 1$.)

This relation correspond to a geometric factor of $A_M(l,m)$.
  
We will describe  how this relation  of $m$, $l$ (near the complete structure)
together  with the Dehn-filling relation
\begin{equation}\label{eq:rel2}
p \log(\frac{m}{m^0}) + q \log(\frac{l}{l^0}) = \frac{\sqrt{-1} \theta}{2}
\end{equation}
gives us Taylor coefficients of $m$ and $l$ 
{\it For simplicity we assume that $m^0 = l^0  = -1$}.

Other cases can be treated in the same way. 

Recall that $a_1$ is not a real number since in
general $a_1$ is the modulus of the flat structure of the
cusp with respect to  $\{\cM ,\cL\}$ when the hyperbolic structure is complete.
 
Since $a_1 \ne -p/q$, as $a_1$ is not real, 
the sets defined by equations \ref{eq:rel1} and \ref{eq:rel2} 
are not tangent at $m^0 = l^0 = -1$. 
We will regard $m$ and $l$ as holomorphic functions
of the variable $\theta$ at $\theta = 0$ locally. 

Our purpose here is to
obtain Taylor coefficients of $m$ and $l$ 
in terms of $\theta$ up to order 3 terms.
We can do so by successive differentiation
of (\ref{eq:rel1}) and (\ref{eq:rel2}) and evaluating at $\theta = 0$.

First we differentiate (\ref{eq:rel2}) to obtain
\begin{equation} \label{eq:rel3}
\frac{p}{m} \, \frac{dm}{d\theta} + \frac{q}{l} \,
\frac{dl}{d\theta} = \frac{\sqrt{-1}}{2}.
\end{equation}
If we evaluate at $\theta = 0 $, we obtain 
\begin{equation} 
\left. p \frac{dm}{d\theta}\right|_{\theta = 0}
+ \left. q \frac{dl}{d\theta}\right|_{\theta = 0} = - \frac{\sqrt{-1}}{2}.  
\end{equation}

On the other hand if we differentiate (\ref{eq:rel1}) and evaluate
at $\theta=0$ , we obtain
\begin{equation}\label{eq:rel4}
\left. \frac{dl}{d\theta}\right|_{\theta = 0} 
= \left. a_1 \frac{dm}{d\theta}\right|_{\theta = 0}.
\end{equation}
From (\ref{eq:rel3}) and (\ref{eq:rel4}) we obtain
$$\frac{dm}{d\theta} \left|_{\theta = 0} = -\frac{\sqrt{-1}}{2(p + a_1 q)},\right. \quad
\frac{dl}{d\theta} \left|_{\theta = 0} 
= -\frac{a_1 \sqrt{-1}}{2(p + a_1 q)} \right.$$ 

Continuing in this manner, we obtain
\begin{equation*}
\begin{split}
&\left. \frac{d^2m}{d\theta ^2} \right|_{\theta = 0}
= \frac{p + (a_1 ^2 + a_2)q}{4(p + a_1 q)^3}, \quad
\left. \frac{d^2 l}{ d \theta ^2} \right|_{\theta = 0}
= \frac{(a_1 - a_2) p + a_1 ^3 q}{4(p + a_1 q)^3},
\\
& \left. \frac{d^3m}{d\theta ^3} \right|_{\theta = 0} \\
=& \frac{\sqrt{-1}\{p^2 +
 (6a_1 ^2 - 2 a_1 ^3 + 6a_2 - 2 a_1 - 3 a_1 a_2 - a_3)pq 
+ (a_1 ^4 + 3 a_1 ^2 a_2 + 3 a_2 ^2 - a_1 a_3) q^2\}}
{8(p + a_1 q)^5},\\
&\left. \frac{d^3l}{d\theta ^3} \right|_{\theta = 0} \\
=& \frac{\sqrt{-1}\{(a_1 - 3 a_2 + a_3)p^2
+ (6a_1 ^3 - 2 a_1 ^4 - 2 a_1 ^2 - 6 a_1 ^2 a_2 - 3 a_2 ^2 + 3 a_1 a_2 + a_1 a_3)pq 
+ a_1 ^5 q^2\}}{8(p + a_1 q)^5}.
\end{split}
\end{equation*}

Recall that $a_2 = a_1 - a_1 ^2$ if the curve represented by  (\ref{eq:rel1}) 
is invariant under the involution $(l,m) \mapsto (1/l , 1/m)$ near 
$(l^0 , m^0)=(-1,-1)$.
Thus we have the following formula
for $m$ , $l$, and $ r \log(\frac{m}{m^0}) + s \log (\frac{l}{l^0})$  
in terms of $\theta$ when $a_2 = a_1 - a_1 ^2$.
\begin{equation}
m = -1 - \frac{\sqrt{-1}}{2(p + a_1 q)} \theta  + \frac{1}{8(p + a_1 q)^2}\theta ^2
+ \sqrt{-1}\frac{p + (3 a_1 - 3 a_1 ^2 + a_1 ^3 - a_3 ) q}{48 (p + a_1 q)^4 } \theta ^3
\label{eq:mlintheta1}
\end{equation}
\begin{equation}
l= -1 - \frac{a_1 \sqrt{-1}}{2(p + a_1 q)} \theta  + \frac{a_1 ^2}{8(p + a_1 q)^2}\theta ^2 +
\sqrt{-1}\frac{(-2 a_1 + 3 a_1 ^2 + a_3 )p + a_1 ^4 q}{48 (p + a_1 q)^4 } \theta ^3
\label{eq:mlintheta2}
\end{equation}
\begin{equation}
r \log(-m) + s \log (-l) = \frac{ \sqrt{-1} (r + a_1 s)}{2 (p + a_1 q)} \theta
+ \frac{\sqrt{-1}(2 a_1 - 3 a_1 ^2 + a_1 ^3 - a_3)(p s - q r)}{48(p + a_1 q)^4} \theta^3
\label{eq:mlintheta3}
\end{equation}
up to order 3 terms.

\subsection{Convergence of the terms of 
Taylor series of geometric factors of A-polynomials}
\label{subsec:convA} 

Let $M$ be a 3-manifold which admits a
double-cusped complete hyperbolic structure and let $\cM_1$ ,
$\cL_1$, $\cM_2 $, and $\cL _2 $ be fixed meridians and longitudes
for the ends of $M$. 
We assume that $M$ can be obtained
from an ideally triangulated hyperbolic manifold by Dehn-filling
some of the ends. Let $\nu$ be the number of the tetrahedra.

We have a holomorphic embedding of an open set $V \subset \bC^2$ 
to $\bC^\nu $ whose image is a subset of $ \cV \subset \bC^\nu$
consisting of points $(z_1 , \cdots, z_\nu )$ satisfying the
gluing consistency relations.(See \cite{NZ} for example). 

Thurston showed that holonomy representations near that of 
complete hyperbolic structure on $M$ have lifts $\rho_0$ to
$\slc$. 
Let $m_1^0, l_1^0 ,m_2^0, l_2^0 $ be the eigenvalues of 
$\rho_0(\cM_1),\rho_0(\cL _1),\rho_0(\cM_2),\rho_0(\cL_2)$. 
(Each of $m_1^0, l_1 ^0 ,m_2^0, l_2^0 $ is either $1$ or $-1$.)
We have a holomorphic map from $\cV$ to $\bC^4$ which assigns to 
each point $\mathbf{z}$ of $\cV$ eigenvalues $ m_1 , l_1 , m_2 , l_2$
of $\rho(\cM_1) , \rho(\cL_1) , \rho(\cM _2) , \rho(\cL_2)$,
respectively, where $\rho$ is a lift of the holonomy representation
corresponding to $\mathbf{z}$. Moreover we choose the
holomorphic map so that the value of $(m_1 , l_1 , m_2 , l_2)$
equals $(m_1^0, l_1^0 ,m_2^0, l_2^0)$ at the point of $\cV$
corresponding to the complete structure.

When the
first end remains a cusp the eigenvalues $m_2 , l_2$
satisfy a relation of the form
\begin{equation*}
\begin{split}
l_2 = l_2 ^0 + a_1 (m_2 - m_2 ^0) 
+ \frac{a_2}{2} (m_2 - m_2 ^0)^2+ \frac{a_3}{6} (m_2 - m_2 ^0 )^3 \\
+ \mbox{higher order terms}.
\end{split}
\end{equation*}

Similarly, when $p_1, q_1$ are coprime integers and $|p_1| + |q_1|$ is large, 
$m_2$ and  $l_2$ satisfy a relation
\begin{equation*}
\begin{split}
l_2 = l_2 ^0 + a_1 ^{p_1 , q_1} (m_2 - m_2 ^0) 
+ \frac{a_2 ^{p_1 , q_1}}{2} (m_2 - m_2 ^0)^2
+ \frac{a_3 ^{p_1 , q_1}}{6} (m_2 - m_2 ^0 )^3 \\
+ \mbox{higher order terms}.
\end{split}
\end{equation*}
when the first end is Dehn-filled along the slope $(p_1 , q_1)$.

\begin{thm} \label{thm:convcoeff}
Let $M$ be a hyperbolic manifold with two cusps. 
Let $M(p_1, q_1)$ be the Dehn-filled $3$-manifold obtained from $M$ 
by a $(p_1, q_1)$-surgery on the first cusp. 
Let $a_i$ be the $i$-th Taylor coefficient of $l_2$ in terms of $m_2$ in 
a geometric factor of the $A$-polynomial of the second cusp of $M$.
Let $a_i^{p_1, q_1}$ be that of $M(p_1, q_1)$.
Then $a_i^{p_1 , q_1} \rightarrow a_i$ 
as $|p_1| + |q_1| \rightarrow \infty$ for $i = 1,2,3$.
\end{thm}

This result shows the convergence of a sequence of coefficients of 
Taylor series of certain geometric components
(defined in Subsection \ref{subsec:A-polynomial}) of the A-polynomials
of the manifolds $M(p_1 , q_1 )$ which are obtained by Dehn
filling a double-cusped hyperbolic manifold $M$ on the first end.
Though we show the convergence of the coefficients up to order
3, we can easily extend our proof to 
show the convergence seems to hold for any higher order terms.

We will prove this theorem in Section \ref{sec:convAproof}.

\section{Maximal tubes in hyperbolic cone
manifolds}\label{sec:maxtubes} 

In this section, we will discuss maximal tubes in hyperbolic cone-manifolds 
in general. We first define the tie classes, the shortest path connecting 
the singularity not homotopic into the singularity. 
We discuss the stability of tie classes under geometric convergence. 
We obtain a formula of the radius of the maximal tube using the traces of 
some elements including the commutator of the tie classes. 
Finally, we obtain the meridian length and the normalized meridian length 
in terms of radius of the maximal tube and cone angles and translation length. 

A hyperbolic manifold of finite volume with a distinguished cusp 
has a horoball neighborhood of the cusp.
The largest of such neighborhoods is called the {\em maximal horoball
neighborhood} of the cusp or the {\em maximal cusp neighborhood}.

Analogously, a hyperbolic cone-manifold whose singular locus is a
knot has standard tube neighborhoods of the singular locus and the
largest of such neighborhoods is called the {\em maximal tube
neighborhood} or the {\em maximal tube around the singular locus}.

Let $M$ be a hyperbolic manifold of finite volume with a cusp. Let
$P$ be a peripheral subgroup of $\pi _1 (M)$. The
boundary of the maximal horoball neighborhood is a torus, say $T$,
which is tangent to itself at a finite number of points. For
each point $x$ of self-tangency, we have a unique geodesic line
which is orthogonal to $T$ at $x$ and tends to the cusp end in
both directions. 

Such a geodesic line corresponds uniquely to an equivalence class 
of the double coset space $P \backslash \pi_1(M)/P$,
i.e., an equivalence class with respect to the equivalence relation $\sim$
defined on $\pi _1(M)$ by $\alpha \sim \beta$ if and only if
$\alpha = \gamma_1 \beta \gamma_2$ or $\alpha = \gamma_1 \beta^{-1} \gamma_2$ 
for some $\gamma_1 , \gamma_2 \in P$.
The class is said to be a {\em tie class} of $P$ or the cusp corresponding 
to $P$. 

Let $\rho _0 : \pi _1 (M) \rightarrow \pslc $ be the holonomy
representation for the hyperbolic structure on $M$ such that 
$\rho_0 (P)$ fixes $\infty$ in the upper space model for the hyperbolic
space $\hs$. Let $H$ be the horizontal plane in $\hs$ which is a
lift of $T$ and let $\tilde{x} \in H$ be a lift of $x$. Then we
have another horosphere $H'$ which is a lift of $T$ and contains
$\tilde{x}$. Then an element $\alpha$ of $\pi_1(M)$ such that
$\alpha (H') = H$ represents the tie class. 
\begin{figure}
\begin{center}\epsfig{file=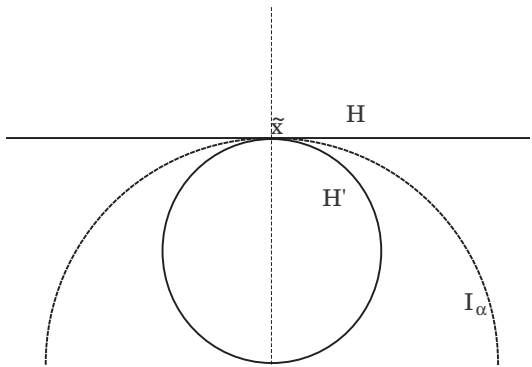, height = 5cm}\end{center}
\caption{A representative of a tie class for maximal horoball
neighborhood} \label{fig:tieforcusp}
\end{figure}

We also characterize the tie class as follows:
A representative of the class has largest isometric
spheres with respect to representations $\rho_0$ with
$\rho_0(P)$ fixing $\infty$.

We define the tie class for
a hyperbolic cone-manifold. Let $N$ be a hyperbolic cone-manifold
whose singular locus $\Sigma$ is a simple closed curve and the
cone angle is less than $2 \pi$. Let $P$ be the peripheral
subgroup of $\pi_1 (N - \Sigma)$. The boundary of the maximal tube
around the singular locus is a torus $T_\Sigma$ with several
points of self-tangency. For each point $x$ of self-tangency, we
have a unique geodesic arc which is the union of two shortest
paths from $x$ to $\Sigma$. Such a geodesic arc corresponds to an
equivalence class in the double coset of $P \backslash \pi _1(N-\Sigma)/P$.

We describe the tie class in the universal covering space: 
Let $x$ be a point of self-tangency on
$T_\Sigma$ and let $\tau$ be one of the shortest paths from $x$ to
$\Sigma$. Let $U_\Sigma$ be the interior of the maximal tube
neighborhood around $\Sigma$ in $N$ minus $\Sigma$  and let
$\widetilde{U}_0$ be the component of the lift of $U_\Sigma$ in
$\widetilde{N - \Sigma}$ which is left invariant by the action of
$P$ on $\widetilde{N - \Sigma}$. Let $\tilde{x _0}$ be a lift of
$x$ in $\widetilde{N - \Sigma}$ lying on the boundary of
$\widetilde{U}_0$ and let $\tilde{\tau_0}$ be the lift of $\tau-\Sigma$ 
in $\overline{\widetilde{U} _0}$ with one end at
$\tilde{x _0}$. By extending the geodesic arc $\tilde{\tau_0}$
past $\tilde{x_0}$ we obtain an open geodesic arc tending to ends of
$\widetilde{N - \Sigma}$ in both directions. One end comes from
$\widetilde{U}_0$ and the other end comes from an image
$\widetilde{U}_0 '$ of $\widetilde{U}_0$ under the action of
$\pi_1 (N-\Sigma)$ on $\widetilde{N - \Sigma}$. Then the tie class
for the maximal tube at $x$ is represented by an element $\alpha$
of $\pi_1 (N-\Sigma)$ such that $\alpha(\widetilde{U}_0 ')=
\widetilde{U}_0$.


Let $\{M_\theta: 0<\theta<\theta_0\}$ be a family of hyperbolic cone-manifolds of
the same topological type such that $M_\theta$ converges (in the
sense of Gromov-Hausdorff) to a complete hyperbolic manifold $M_0$
of finite volume with at least one cusp and, for each $\theta$,
$M_\theta$ has a singular locus $\Sigma_\theta$ with cone angle
$\theta$ and $M_\theta - \Sigma_\theta$ is homeomorphic to a fixed $3$-manifold 
$M_0$ by a map $\phi _\theta: M_0 \rightarrow M_\theta -\Sigma_\theta$. 
Then we have holonomy representations $\rho_0:\pi_1 (M_0) \rightarrow \pslc$ and
$\rho_\theta :\pi_1(M_\theta - \Sigma_\theta) \rightarrow \pslc$
for $0<\theta<\theta_0$ so that 
$\rho_\theta \circ (\phi_\theta)_* \rightarrow \rho_0$ is an isomorphism. 
Let $P$ be a peripheral subgroup
of $\pi_1(M_0)$. 

\begin{prop}\label{prop:tieclass}
Assume the notations above.
When $M_0$ has a tie class, say
$[\alpha]$ with $\alpha \in \pi_1(M_0)$, for the cusp corresponding to $P$, 
the cone manifold
$M_\theta$ also has a tie class with respect to 
$(\phi_\theta)_*(P)$ which is represented by $(\phi _\theta)_* (\alpha)$
when $\theta$ is small. 
\end{prop}
\begin{proof}
This is obvious from Gromov-Hausdorff topology since a sequence of maximal 
tubes must converge to a maximal horoball.  
\end{proof}


\begin{prop}\label{prop:formularadius}
Let $N$ be a hyperbolic cone manifold whose singular locus is a
simple closed curve and the cone angle is less than $2\pi$. Let
$P$ be a peripheral subgroup of $\pi_1 (N-\Sigma)$ and $\rho:
\pi_1(N-\Sigma) \rightarrow \pslc$ be a  holonomy representation
of the hyperbolic structure on $N-\Sigma$. Let $\alpha \in
\pi_1(N-\Sigma)$ represent a tie class for the maximal tube. Then
the radius $R$ of the maximal tube around the singular locus
satisfies
$$\cosh(2R) =\frac{|\tr\rho(\alpha \gamma \alpha ^{-1} \gamma^{-1}) - 2|
+|\tr ^2 \rho(\gamma) - \tr\rho(\alpha \gamma \alpha ^{-1}
\gamma^{-1}) - 2| }{|\tr ^2 \rho(\gamma)-4 |},$$ where $\gamma$ is
any element of $P$ such that $\rho(\gamma) \ne I$.
\end{prop}
\begin{proof}

Let ${\dev}$ be the developing map $\widetilde{N - \Sigma}
\rightarrow \hs$ for the hyperbolic structure on $N-\Sigma$ such
that ${\dev}\circ \gamma = \rho(\gamma) \circ {\dev}$ for
any $\gamma \in \pi_1(N - \Sigma)$. Let $x$, $\tilde{x _0}$,
$\tau$, $\tilde{\tau_0}$, $\bU_\Sigma$, and $\widetilde{U} _0$ be
as above. Then $\rho(P)$ fixes a geodesic line $\tilde{\Sigma} $
in $\hs$ and ${\dev}(\widetilde{U} _0)$ is the set of points in
$\hs$ lying within the radius of the maximal tube from
$\tilde{\Sigma} $. Let $\tilde{x}$ and $\tilde{\tau}$ be the
images of $\tilde{x _0}$ and $\tilde{\tau _0}$ under 
$\dev|{\overline{\tilde {U_0}}}$, respectively. Then
$\dev(\widetilde{U} _0)$ and ${\dev}(\alpha^{-1}(\widetilde{U}_0))$ 
are hyperspherical regions around $\tilde{\Sigma} $ and
$\rho(\alpha ^{-1})(\tilde{\Sigma})$, respectively which are
tangent at $\tilde{x}$.

\begin{figure}
\begin{center}
\epsfig{file=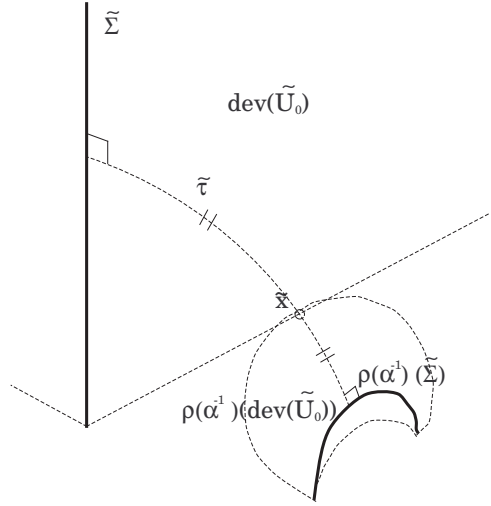, height = 7cm}
\end{center}
\caption{A representative of a tie class for maximal tube}
\label{fig:maxtube}
\end{figure}

Recall that $R$ is half the
distance between $\tilde{\Sigma}$ and
$\alpha^{-1}(\tilde{\Sigma})$, and $\tilde{\Sigma} $ is the axis
of $\rho(\gamma)$ for any $\gamma \in P$ such that $\rho(\gamma)
\ne I$.  
In the upper half space model, geodesic lines are
represented by
pairs of different extended complex numbers. We may assume that $\tilde{\Sigma} $
is the vertical geodesic represented by $(0, \infty)$ and that
$\rho(\gamma)$ is of the form
$$ \begin{bmatrix} u & 0 \\ 0 & 1/u \end{bmatrix} .$$
Let 
$$\rho(\alpha)= \begin{bmatrix} a & b \\ c & d
\end{bmatrix}$$ 
with $ad - bc = 1$.
Then the geodesic line
$\rho (\alpha^{-1} )(\tilde{\Sigma})$ is represented by $(-b/a ,-d/c)$. 
So $\cosh$ of the distance between $\tilde{\Sigma}$ and
$\rho (\alpha^{-1})(\tilde{\Sigma})$ is $|ad| + |bc| = |bc| + |bc+ 1|$ 
by the following lemma.
Since $\tr\rho(\alpha \gamma \alpha ^{-1} \gamma^{-1}) -2 
= -bc ( u - 1/u)^2 = -bc (\tr ^2 \rho(\gamma)-4)$. 
Now it is straightforward
to check the equality in the proposition. 
\end{proof}

\begin{lem}
The distance $d$ between two geodesic lines $(w_1, w_2)$ and
$(w_3, w_4)$ satisfies
$$\cosh {d} = \frac {1 + |[w_1,w_2;w_3,w_4]|}{|1-[w_1, w_2; w_3, w_4]|} \; , $$ 
where $$[w_1,w_2;w_3,w_4] = \frac{(w_1 - w_3)(w_2 - w_4)}{(w_1 - w_4)(w_2 - w_3)}$$ 
is the cross-ratio.
\end{lem}
\begin{proof}
Since the cross-ratios and hyperbolic distances are
invariant under hyperbolic isometries, we need only to prove the
lemma when $w_1 = -1 , w_2 = 1 , w_3 = -w,  w_4 = w$ for some
$w \in \bC$. But in this case the equality is easily shown.
\end{proof}


Let $\theta$ be the cone angle of the singular locus and let $\gamma_0$ 
be an element of $P$ such that $\rho(\gamma_0)$ is not
elliptic and moves by minimal distance along its axis. Let
$t$ be the axis length of $\rho(\gamma_0)$ which equals the
absolute value of the real part of $2 \log(u)$ when
$\tr\rho(\gamma_0) = u + 1/u$ , i.e., the length
of the singular locus. Let $\mu$ and $\hatmu$ be the length and
the normalized length $\hatmu$ of the meridian curve on the
boundary of the maximal tube, respectively. 
$$\mu = \theta
\sinh (R)= \theta \sqrt{\frac{\cosh(2R) - 1}{2}}$$
and
\begin{equation}\label{eq:hatmu}
\hatmu^2=\frac{\theta \tanh(R)}{t} = \frac{\theta}{t}
\sqrt{\frac{\cosh(2R) - 1}{\cosh(2R) +1}},
\end{equation}
we can
express $\mu$ and $\hatmu$ in terms of the cone angle and the
traces of holonomy images of certain elements of $\pi_1(N -\Sigma)$ 
related to the tie classes.


\section{Maximal Tubes in Whitehead Link Cone Manifolds}\label{sec:main}

The purpose of this section is prove Theorem \ref{thm:main}. 
The outline of the proof will be given after the theorem is stated.

\begin{figure}
\begin{center}
\epsfig{file=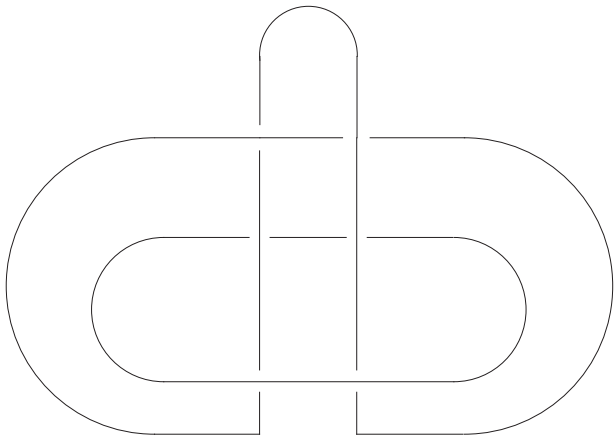, height=6cm}\end{center} \caption{The Whitehead Link} 
\label{fig:wh1}
\end{figure}

We choose the meridian-longitude pair for each
end. For each coprime pair $p_1$, $q_1$ of integers, let $W(p_1, q_1)$ 
be the manifold obtained from the Whitehead link complement $W$ 
by Dehn-filling the first
end along the slope $(p_1 , q_1 )$. If $|p_1| + |q_1|$ is
sufficiently large, $W(p_1 , q_1)$ is hyperbolic, and $\theta > 0$ is
small, let $W(p_1 , q_1)(p_2 , q_2; \theta)$ be the hyperbolic
cone manifold obtained by generalized Dehn-filling of 
$W(p_1 , q_1)$ on the second end along the slope $(p_2,q_2)$ with cone
angle $\theta$ and let $\mu = \mu_{p_1, q_1, p_2, q_2} (\theta)$ be
the length of the surgery curve on the maximal tube of 
$W(p_1 , q_1)(p_2 , q_2; \theta)$ around the singularity. 
Let $\hat{\mu}=\hat{\mu}_{p_1, q_1, p_2, q_2} (\theta)$ be the normalized length
of the surgery curve, that is, $\mu$ divided by the square root of the area of the
boundary of the maximal tube. 

\begin{thm} \label{thm:mainthm}
Following the above notations, 
if $|p_1| + |q_1|$ is sufficiently large, then for any coprime
pair $p_2,q_2$ of integers except for at most one, $\hat{\mu}$ for 
$W(p_1, q_1)(p_2, q_2; \theta)$ is decreasing and
$\hat{\mu}^2 + \theta ^2 $ is increasing near $\theta=0$.
\end{thm}

Let us give an outline of the proof of Theorem \ref{thm:mainthm}.
We start with basic materials about the deformation of hyperbolic
structures on the Whitehead link complement $W$ near the complete
structure: 

In Subsection \ref{subsec:NZforW} we present a
parametrization for the hyperbolic structures near the complete
structure using a decomposition of $W$ into ideal tetrahedra and
in Subsection \ref{subsec:holW} we obtain all the lifts of holonomy
representations (up to conjugacy) for hyperbolic structures on
$W$ near the complete structure. 

In Subsection \ref{subsec:maxtubesW},
we obtain formulas for $\hat{\mu}_{p_1, q_1, p_2, q_2}^2$ in terms of the
eigenvalues $m_2$ and $l_2$ of
of holonomy (for the hyperbolic structure on $W$ inducing 
the hyperbolic cone-structure on
$W(p_1 , q_1)(p_2 , q_2; \theta)$) images of $\cM_2$ and
$\cL_2$ using the results of Subsection \ref{subsec:holW}. 

The needed  observation is that each $W(p_1 , q_1)$ has a unique tie class
for maximal cusp which comes from the tie class for the maximal
cusp for the second end in $W$ and this tie class is easily computable.

In our case, 
if $M$ has a unique tie class for its maximal cusp, we can
find a representative for the tie class, and 
we can compute a geometric component of the $A$-polynomial zero-set,
then we can compute the normalized length of the meridian curve
in $M(p,q;\theta)$ in terms of $\theta$. 

As explained above, we know that
the single-cusped complete hyperbolic manifold
$W(p_1 , q_1)$(when $|p_1| + |q_1|$ is large), 
has a unique tie class for its maximal cusp and we can
easily find a representative for the tie class. 

In Subsection \ref{subsec:completion}, putting all
things together, we will obtain Taylor coefficients 
of the function $\hat{\mu}_{p_1, q_1, p_2, q_2}^2 (\theta)$ 
up to order 2 terms. 
Using concrete values of $a_1 , a_2 ,a_3$ we
can show Theorem \ref{thm:mainthm} 
when the first cusp remains a cusp, i.e., $p_1 = q_1 =\infty$. 
Then we will complete the proof using the convergence 
$a_i ^{p_1 , q_1} \rightarrow a_i$ $(i=1,2,3)$ proved 
in the last section \ref{sec:convAproof}.

\subsubsection{The parametrization of hyperbolic
structures on the Whitehead link complement} \label{subsec:NZforW}

\begin{figure}
\begin{center}
\epsfig{file=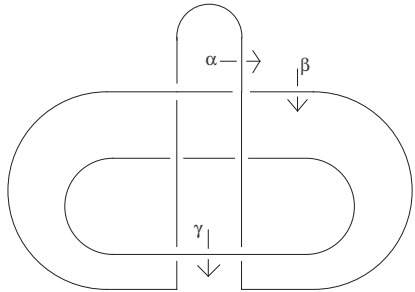, height=4cm}
\end{center}
\caption{Whitehead link complement with Wirtinger generators}
\label{fig:wh2}
\end{figure}

We start with a holonomy representation and a  picture of the
Ford domain for the complete hyperbolic structure on Whitehead
link complement. We use the presentation
$$\langle \alpha,\beta,\gamma:\alpha \gamma = \gamma \beta , \,
\gamma \alpha \beta \alpha ^{-1}= \alpha \beta^{-1} \alpha \beta
\alpha ^{-1} \gamma \rangle$$
of the fundamental group of $W$ coming from Figure \ref{fig:wh2}.
We fix the following system of meridians and longitudes for ends
of $W$: $\cM_1 = \gamma , \cL_1 = \alpha \beta^{-1} \alpha^{-1}
\beta , \cM_2 = \alpha , \cL_2 = \gamma \alpha^{-1} \gamma^{-1}
\alpha \beta^{-1} \alpha$.

The group homomorphism $\rho_0 : \pi _1 (W) \rightarrow \slc $
which sends $\alpha , \beta , \gamma$ onto 
$$
\begin{bmatrix}-1&1\\0&-1 \end{bmatrix},
\quad \begin{bmatrix}-1&0\\2\sqrt{-1}&-1 \end{bmatrix},
\quad\begin{bmatrix}-2&-(1+\sqrt{-1})/2\\1-\sqrt{-1}&0
\end{bmatrix},
$$ 
respectively, is the lift of a holonomy
representation of the complete hyperbolic structure on $W$. 

\begin{figure}
\begin{center}
\epsfig{file=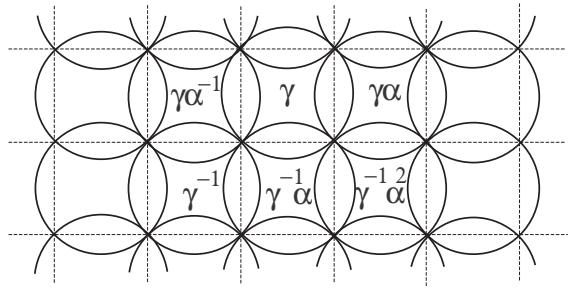, height=8cm,angle=90}
\end{center}
\caption{Ford domain for the second cusp of Whitehead link
complement} \label{fig:whfd4}
\end{figure}

The Ford domain with respect to the second cusp corresponding to
$\rho_0$ together with isometric spheres is described in Figure
\ref{fig:whfd4}.
From the description of the Ford domain for the second end, we
obtain a decomposition of W into four ideal tetrahedra.
\begin{figure}
\begin{center}
\epsfig{file=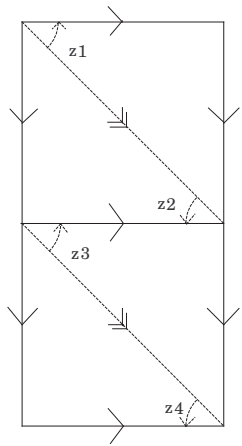, height=6cm}
\end{center}
\caption{A decomposition of $W$ into ideal tetrahedra}
\label{fig:idealdecomp}
\end{figure}

\begin{figure}
\begin{center}
\epsfig{file=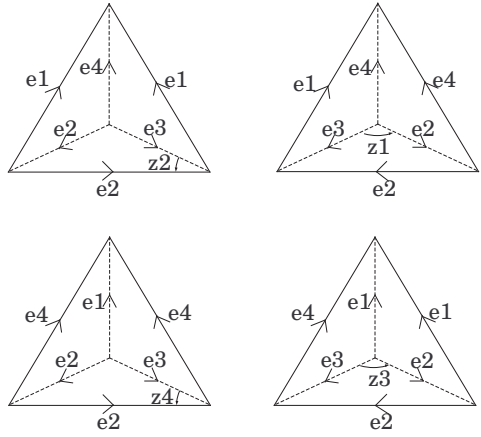, height=6cm}
\end{center}
\caption{Gluing pattern for the tetrahedra}
\label{fig:triangulation}
\end{figure}

Now let $z_1, z_2, z_3, z_4$ be the parameters corresponding to
the four tetrahedra as described in Figure
\ref{fig:triangulation}. 
$z_1 = z_2 = z_3 = z_4 = (1+\sqrt{-1})/2$ 
corresponds to the complete hyperbolic structure.
For the induced metric to be nonsingular along the edges
$e_1,e_2,e_3,e_4$, $z_i$ 's must satisfy the following gluing
consistency relations.
\begin{figure}
\begin{center}
\epsfig{file=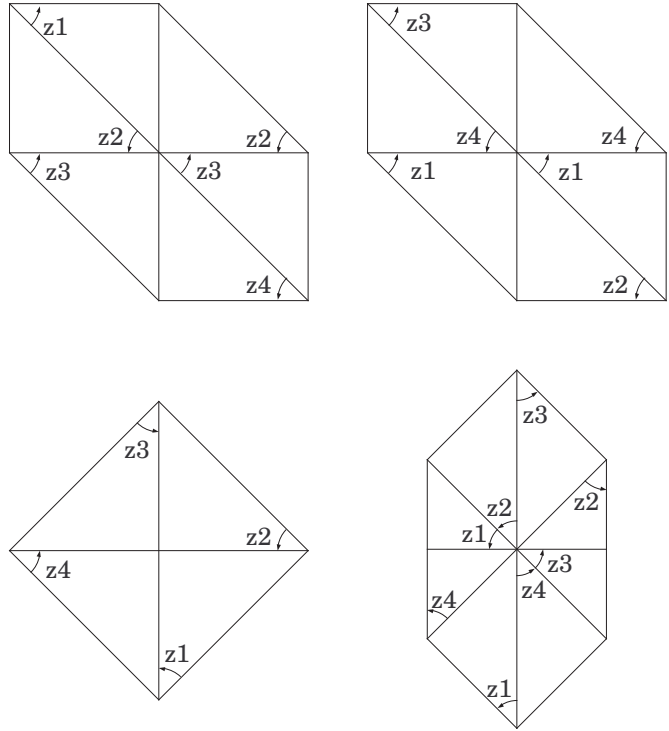, height=10cm}
\end{center}
\caption{Gluing consistency for edges of the triangulation}
\label{fig:consistency}
\end{figure}
\begin{eqnarray} \label{eq:consistencyforW1}
(1-z_1)(1-z_4) = (1-z_2)(1-z_3)\phantom{aaa}\;\\
(1-z_1)(1-z_2)(1-z_3)(1-z_4)= z_1 z_2 z_3 z_4
\label{eq:consistencyforW2}
\end{eqnarray}
For $z_i$'s (all near $z^0$) satisfying those relations,  we have
the following relation between $z_i$'s and the eigenvalues of the
corresponding $\slc$-holonomy images of the meridian and longitude
for the second end.
\begin{figure}
\begin{center}
\epsfig{file=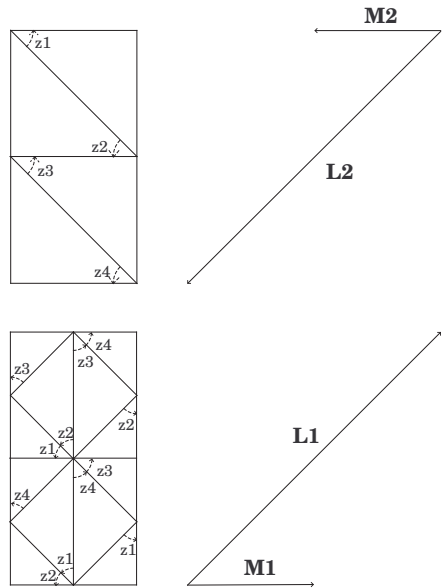, height=8cm}
\end{center}
\caption{Holonomy images of meridians and longitudes}
\label{fig:holonomyml}
\end{figure}
\begin{eqnarray}\label{eq:mlinz1}
m_1 = -\sqrt{\frac{1-z_4}{1-z_2}} \quad \left(\;= -\sqrt{\frac{1-z_3}{1-z_1}} \; \right),\\
l_1 = -\frac{1-z_4}{1-z_2}\sqrt{\frac{z_3 z_4}{z_1 z_2}} 
\quad \left( \;= -\frac{1-z_3}{1-z_1}\sqrt{\frac{z_3 z_4}{z_1 z_2}}\; \right),\\
m_2 = -\sqrt{\frac{1-z_2}{1-z_1}}\quad \left(\;= -\sqrt{\frac{1-z_4}{1-z_3}}\; \right), \\
l_2 = -\frac{1-z_2}{1-z_1}\sqrt{\frac{z_2 z_4}{z_1 z_3}} \quad
\left( \;= -\frac{1-z_4}{1-z_3}\sqrt{\frac{z_2 z_4}{z_1 z_3}}\; \right),  
\label{eq:mlinz4}
\end{eqnarray}
where the branch of square root function is chosen to take value 1 at 1.

\subsubsection{Holonomy representations of the Whitehead link
group} \label{subsec:holW} 

We claim that for any $x, y \in \bC$, there
is a homomorphism $\rho: \pi _1(W) \rightarrow \slc$ such that
\begin{eqnarray*} \rho(\alpha) =\begin{bmatrix}x & 1 \\
0 & 1/x\end{bmatrix},\phantom{aaaaaaa}
\rho(\beta)= \begin{bmatrix}x & 0 \\ y & 1/x\end{bmatrix},&\\
\rho(\gamma)= 
\begin{bmatrix}-\frac{x\{x^2 y^2 + x^2 (x^2 - 3) y - (x^2 - 1 )^2\}}{z w} & z \\
z y & \frac{z(1-x^2)}{x} \end{bmatrix},& \end{eqnarray*} where 
$w = (x^2 y - x^2 + 1)\{x^2 y + (x^2 - 1)^2\},\quad z=
\sqrt{x^2(1-x^2-y)/w}$ and that the lift of any holonomy
representation for $\pi_1(W)$ near $\rho _0$ is conjugate to one
of this form (with $(x,y,z)$ near $(-1, 2 \sqrt{-1} , -( 1 + \sqrt{-1})/2 )$).\\

\begin{proof}[Proof of the claim] 
Let $\rho : \pi_1 (W) \rightarrow \slc$ be a
lift of a  holonomy representation such that $\rho(\alpha)$ and
$\rho(\beta)$ does not commute. Then by conjugation in $\slc$ we
may assume that $\rho(\alpha)$ and $\rho(\beta)$ are of the form
\begin{equation}\label{eq:alphabeta} 
\begin{bmatrix}x & 1 \\ 0 & 1/x\end{bmatrix},\quad
\begin{bmatrix}x & 0 \\ y & 1/x\end{bmatrix},
\end{equation}
respectively, since $\rho(\alpha)$ and $\rho(\beta)$ are conjugate in
$\slc$. By the relation $\alpha \gamma = \gamma \beta$,  $\rho(\gamma)$ must
be of the form
\begin{equation}\label{eq:gamma}
\begin{bmatrix}
\frac{x(1 + y z^2)}{z(1 - x^2)} & z \\ z y & \frac{z(1-x^2)}{x}
\end{bmatrix}
\end{equation}
Finally from the relation 
$\gamma \alpha \beta \alpha ^{-1}= \alpha \beta^{-1} \alpha \beta
\alpha ^{-1} \gamma$, we obtain $z$ in terms of $x,y$
 as in the claim.
\end{proof}


For a given holonomy representation $\rho$ described above, 
we take $x$ as the eigenvalue $m_2$ of the holonomy image of
the meridian $\alpha= \cM_2$ for the second end. 
Since
$$
\rho(\cL _2) = 
\begin{bmatrix} \frac{-1 + x^2 - x^2 y}{-1+ x^2 + y}
& \ast \\ 0 & \frac{-1+ x^2 + y}{-1 + x^2 - x^2 y}
\end{bmatrix}
$$ 
we should take the eigenvalue $l_2$ of the holonomy image of the longitude $\cL_2$
as
\begin{equation} \label{eq:l2inXY}
l_2 = \frac{-1 + x^2 - x^2 y}{-1+ x^2 + y}.
\end{equation}
Then we obtain the following relation between eigenvalues $m_1,l_1,m_2,l_2$
of the holonomy images 
$\rho(\cM _1),\rho(\cL _1),\rho(\cM_2),\rho(\cL_2)$ by computing 
$\tr^2 \rho(\cM_1)$ and $\tr \rho(\cL_1)$ in terms of $x , y ,z$ in the claim above.
\begin{eqnarray*}
(m_1 + \frac{1}{m_1})^2& =
& \frac{(1+l_2)^2 (m_2 ^4 - l_2)}{l_2(l_2 + m_2 ^2)(m_2 ^2 - 1)} \\
l_1 + \frac{1}{l_1}\phantom{aaa} & 
=& \frac{l_2 ^2(1+ m_2 ^4 ) 
+ l_2 (-1 + 2 m_2 ^2 + 2 m_2 ^4 + 2 m_2 ^6 - m_2 ^8) + m_2 ^4 + m_2 ^8}
{m_2 ^2 ( l_2 + m_2 ^2 )^2}
\end{eqnarray*}

\subsection{Maximal tubes in Whitehead link cone manifolds}
\label{subsec:maxtubesW}

The elements of $\pi _1(W)$ which represent a unique tie class
for the maximal cusp of
$W$ with respect to the peripheral subgroup $P_2$ for the second
end containing $\alpha \in \pi_1(W)$ projects to an element of
$\pi_1(W(p_1,q_1)(p_2,q_2;\theta) -
\Sigma(p_1,q_1)(p_2,q_2;\theta))$ representing a tie class when
$|p_1| + |q_1|$ is large and $\theta$ is small. 

From the
description of Ford domain of $W$ with respect to the second end
given in Subsection \ref{subsec:NZforW} we see $\gamma \in \pi_1(W)$
represents the unique tie class for the second cusp with respect
to $P_2$. Also we see that when $|p_1| + |q_1|$ is large and
$W(p_1,q_1,\infty,\infty)$ is hyperbolic, $\gamma$
still represents the unique tie class for the single-cusped complete
hyperbolic manifold $\pi_1 (W(p_1,q_1,\infty,\infty))$ by considering the isometric
spheres of images of a suitable holonomy representation. 

From results of Section \ref{sec:maxtubes}
(Proposition \ref{prop:formularadius} and (\ref{eq:hatmu})), we see that for
the hyperbolic cone manifold $W(p_1,q_1)(p_2,q_2;\theta)$ the
normalized length $\hatmu$ of the surgery curve on the boundary
of the singular locus satisfies
\begin{equation} \label{eq:hatmuW1}
\begin{split}
\hatmu^2 = \frac{\theta}{2 |\re (r_2 \log(-m_2) + s_2 \log(-l_2))|}
\phantom{aaaaaaaaaaaaaaaaaaaaaaaaaaaaaaaaaa}\\
\times  \sqrt {\frac{|\tr(\rho(\alpha \gamma \alpha^{-1}\gamma^{-1}) - 2)|
+ |\tr ^2 ( \rho(\alpha)) - \tr ( \rho ( \alpha \gamma \alpha^{-1}\gamma^{-1} ) - 2|
- |\tr ^2(\rho(\alpha)) - 4|}
{|\tr(\rho(\alpha \gamma \alpha^{-1}\gamma^{-1}) - 2| + |\tr ^2(\rho(\alpha))
-\tr(\rho(\alpha \gamma \alpha^{-1}\gamma^{-1}) - 2| + |\tr ^2(\rho(\alpha)) -  4|} }
\end{split}
\end{equation}
  where 
$\rho: \pi_1(W(p_1,q_1)(p_2,q_2;\theta) - \Sigma(p_1,q_1)(p_2,q_2;\theta)) \rightarrow \slc$ 
is a lift of a holonomy representation for $W(p_1,q_1)(p_2,q_2;\theta)$ and
$r_2 , s_2$ are integers satisfying $p_2 s_2 - q_2 r_2 =1$,   
when $|p_1|+ |q_1|$ is large and $\theta$ is small.

  We obtain that 
$\tr\rho(\alpha \gamma \alpha^{-1}\gamma^{-1}) = \tr \rho(\alpha \beta^{-1}) = 2 - y$
from equations \ref{eq:alphabeta} and \ref{eq:gamma}.
Using (\ref{eq:l2inXY}), we obtain
\begin{equation}\label{eq:trcomm}
\tr\rho(\alpha \gamma \alpha^{-1}\gamma^{-1}) - 2 
= -\frac{(-1 + m_2 ^2)(1 - l_2)}{l_2 + m_2 ^2 }.
\end{equation}

\subsection{Completion of the proof of the main result}
\label{subsec:completion}

Let $W$ be the Whitehead link complement and let 
$\cM_1 , \cL_1 , \cM_2 , \cL_2 $ 
be the meridians and longitudes for the ends of $W$ as before. 

From (\ref{eq:hatmuW1}) and
(\ref{eq:trcomm}), we have
\begin{equation} \label{eq:hatmuW2}
\begin{split}
\hatmu ^ 2 =& \frac{\theta}{2 |\re (r_2 \log(-m_2) + s_2 \log(-l_2))|}  \\
\times & \sqrt{\frac{|\frac{(m_2 ^2 - 1)(1 - l_2)}{m_2 ^2 + l_2}|
+ |(m_2 - \frac{1}{m_2})^2 +
\frac{(m_2 ^2 - 1)(1 - l_2)}{m_2 ^2 + l_2}|- |(m_2 - \frac{1}{m_2})^2|}
{|\frac{(m_2 ^2 - 1)(1 - l_2)}{m_2 ^2 + l_2}|+ |(m_2 - \frac{1}{m_2})^2
+ \frac{(m_2 ^2 - 1)(1 - l_2)}{m_2 ^2 + l_2}|+ |(m_2 - \frac{1}{m_2})^2|} } \; ,
\end{split}
\end{equation}
where $r_2 , s_2$ are integers satisfying 
$p_2 s_2 - q_2 r_2 =  1$,
when $|p_1|+ |q_1|$ is large and $\theta$ is small.

When $p_1 = q_1 = \infty$, we have a Taylor expansion
$$l_2 = -1 + a_1 (m_2 + 1) + \frac{a_2}{2} (m_2 + 1)^2+ \frac{a_3}{6} (m_2 + 1 )^3 +
O((m_2 + 1)^4),$$
where $a_1 = 2 + 2 i , \; a_2 = 2 - 6 i , \; a_3 = -12 $, from
the $A$-polynomial
$$-l_2 + l_2 ^2 + 4 l_2 m_2 + m_2 ^4 -l_2 m_2 ^4 $$
of the manifold $W(\infty, \infty)$ with respect to 
$\{\cM _2 , \cL_2\}$.
So when $|p_1|+ |q_1|$ is large, we have a Taylor expansion
$$l_2  = -1 + a_1 ^{p_1 , q_1}(m_2 + 1) + \frac{a_2 ^{p_1 , q_1}}{2} (m_2 + 1)^2+
\frac{a_3 ^{p_1 , q_1}}{6} (m_2 + 1 )^3 +   O((m_2 + 1)^4) ,$$
for the manifold $W(p_1 , q_1)$
with $a_i ^{p_1 , q_1} \rightarrow a_i$ as $|p_1| + |q_1| \rightarrow \infty$ for $i=1,2,3$ 
as in Subsection \ref{subsec:convA}.
By Lemma \ref{lem:involution}, 
$a_2 ^{p_1 , q_1} = a_1 ^{p_1 , q_1} - (a_1 ^{p_1 , q_1})^2$ 
when $|p_1|+ |q_1|$ is large.

Using (\ref{eq:mlintheta1})--(\ref{eq:mlintheta3}), we obtain
$$
(m_2 - \frac{1}{m_2})^2 = -\frac{\theta ^2}{(p_2 + a_1 ^{p_1 , q_1} q_2  )^2} + O(\theta ^3)$$
$$\frac{(m_2 ^2 - 1)(1 - l_2)}{m_2 ^2 + l_2} = \frac{4}{2 - a_1 ^{p_1 , q_1}}
+ \frac{-6 a_1 ^{p_1 , q_1} + 3 (a_1 ^{p_1 , q_1} )^2 - 2 a_3 ^{p_1 , q_1}}
{12 (2 - a_1 ^{p_1 , q_1} )^2 (p_2 + a_1 ^{p_1 , q_1} q_2  )^2 }\theta^2  + O(\theta ^3) $$ and
\begin{equation*}
\begin{split}
\re (r_2 \log(-m_2) + s_2 \log(-l_2))
\phantom{= \re (r_2 \log(-m_2) + s_2 \log(-l_2)) }
\phantom{aaaaaaaaaa}
\\= - \frac{\im \; a_1 ^{p_1 , q_1}}{2|p_2 + a_1 ^{p_1 , q_1} q_2|^2}\theta
\phantom{\re (r_2 \log(-m_2) + s_2 \log(-l_2))} \phantom{aaaaaaaaaaaaaaaaa}
\\- \frac{\im\{(2 a_1 ^{p_1 , q_1} - 3 (a_1 ^{p_1 , q_1})^2 + (a_1 ^{p_1 , q_1})^3 - a_3  ^{p_1 , q_1} )
( p_2   + \overline{a_1 ^{p_1 , q_1}} q_2 )^4 \}}{48 |p_2 + a_1 ^{p_1 , q_1} q_2| ^4}\theta^3 + O(\theta ^4 )
\end{split}
\end{equation*}

From (\ref{eq:hatmuW2}), we obtain
\begin{equation*}
\begin{split}
& \hatmu_{p_1,q_1,p_2 , q_2}^2 (\theta) \\
=& \frac{|p_2 + a_1 ^{p_1 , q_1} q_2|^2} {|\im \; a_1 ^{p_1 , q_1} |} +
\frac{B_{0} ^{p_1 , q_1}  p_2 ^4 + B_{1} ^{p_1 , q_1} p_2 ^3 q_2
+  B_{2} ^{p_1 , q_1} p_2 ^2 q_2 ^2 + B_{3} ^{p_1 , q_1} p_2
q_2 ^3 + B_{4} ^{p_1 , q_1} q_2 ^4}{|p_2 + a_1 ^{p_1 , q_1} q_2|^4} \theta ^2
\end{split}
\end{equation*}
up to order 2 terms when $|p_1|+ |q_1|$ is large, where  each $B_j ^{p_1 , q_1} $ is a constant
which depend continuously on $a_1 ^{p_1 , q_1}, a_2 ^{p_1 , q_1}, a_3 ^{p_1 , q_1}$, 
that is, $$B_j ^{p_1 , q_1} = B_j (a_1 ^{p_1 , q_1}, a_2 ^{p_1 , q_1}, a_3 ^{p_1 ,  q_1})$$ 
for a continuous function $B_j$ which is defined on a neighborhood of $(a_1, a_2 , a_3 )$ in $\bC^3$.
In the same way, we have
\begin{equation} \label{eq:firstendcuspidal}
\begin{split}
& \hatmu_{\infty,\infty, p_2 , q_2}^2 (\theta) \\=& \frac{|p_2 + a_1 q_2|^2} {|\im \; a_1 |} +
\frac{B_{0} ^{\infty, \infty}  p_2 ^4 + B_{1} ^{\infty, \infty} p_2 ^3 q_2
+  B_{2} ^{\infty, \infty} p_2 ^2 q_2 ^2 + B_{3} ^{\infty, \infty} p_2
q_2 ^3 + B_{4} ^{\infty, \infty} q_2 ^4}{|p_2 + a_1  q_2|^4} \theta ^2
\end{split}
\end{equation}
up to order 2 terms, where $B_{j} ^{\infty, \infty}= B_j (a_1 , a_2, a_3)$. 
Using the values $a_1, a_2 , a_3$, we obtain
$ B_{0} ^{\infty, \infty} = -1/12, \; 
B_{1} ^{\infty, \infty} = -2/3,\; 
B_{2} ^{\infty, \infty} = -4 ,\; 
B_{3} ^{\infty, \infty} = -32/3, \; 
B_{4} ^{\infty, \infty}$ $= -32/3$.

All the values of the function 
$$\frac{-x^4 - 8 x ^3 - 48 x^2 - 128 x - 128}{12(x^2 + 4 x + 8)^2}(-\infty < x <\infty)$$ 
lie in the interval $[-1/6, -1/12] \subset (-1 , 0)$.
Thus the coefficient of the second order term in
(\ref{eq:firstendcuspidal}),
\begin{equation}\label{eq:mainfo}  
\frac{-p_2 ^4 - 8 p_2 ^3 q_2 - 48p_2 ^2 q_2 ^2 - 128 p_2 q_2 ^3 - 128 q_2 ^4}
{12(p_2 ^2  + 4 p_2 q_2 + 8 q_2 ^2)^2} 
\end{equation}
lies in $[-1/6, -1/12] \subset (-1 , 0)$ for any coprime pair of integers
$p_2 , q_2$. This completes the part of our proof when $(p_1 , q_1)$ 
is the infinity. Since $a_1 ^{p_1 , q_1} \rightarrow a_1$, and 
$B_j ^{p_1 , q_1}\rightarrow  B_j ^{\infty,\infty} $ for $j=0 , \dots, 4$ 
as $|p_1| + |q_1| \rightarrow  \infty$, our proof for $\hatmu^2$
is completed. For $\hatmu^2 + \theta^2$, the proof follows 
since the range of equation \ref{eq:mainfo} plus $1$ lies in 
$[5/6, 11/12]$ which is positive.

\section{Convergence of the terms of Taylor series 
of geometric factors of A-polynomials}
\label{sec:convAproof} 

In this section, we will provide a proof of Theorem \ref{thm:convcoeff}.
We assume the notations of Subsection \ref{subsec:convA}.

Let $M$ be a complete hyperbolic manifold with two cusps and 
decomposes into $\nu$ ideal tetrahedra. 
Let $\mathbf{z}^0= (z_1 ^0, \cdots, z_\nu ^0) \in \bC^\nu$ correspond to the
complete structure. Let $\cV$ be the intersection of a small
neighborhood of $\mathbf{z}^0$ in $\bC^\nu$ and the set of points
satisfying the gluing consistency relations as above . 

By results
of Neumann and Zagier\cite{NZ}, there is a holomorphic embedding
$\iota$ of a small open subset $V$ of $\bC^2$ onto $\cV$ (by
shrinking $\cV$ if necessary). 

Let $\rho_0: \pi_1 (M) \rightarrow \slc$ be 
the lift a holonomy representation of the complete
structure and let $m_1 ^0, l_1 ^0 ,m_2 ^0, l_2 ^0 $ be the
eigenvalues of $\rho(\cM _1),\rho(\cL _1),\rho(\cM_2),\rho(\cL_2)$. 

As mentioned before, by suitable choices of
the eigenvalues $m_1, l_1, m_2, l_2$ of 
$\rho(\cM_1), \rho(\cL_1),\rho(\cM_2),\rho(\cL_2)$ 
for the lifts $\rho$ of holonomy representations for
nearby hyperbolic structures  we obtain a holomorphic map $G: \cV\rightarrow \bC^4$ 
which can be regarded as a pair of maps $G_1, G_2:\cV \rightarrow \bC^2$.
$G_1 = (m_1, l_1)$ and $G_2 = (m_2, l_2)$. 

Moreover we can choose $G$ so that
$G_1(\mathbf{z}^0) = (m_1^0 , l_1^0)$ and $ G_2(\mathbf{z}^0) =(m_2^0 , l_2^0)$. 
Let $\cC_1 \subset \cV$ be $G_1 ^{-1}\{(m_1 ^0, l_1 ^0)\} \cap \cV$ and $C_1 \subset V$ 
be $\iota^{-1}(\cC_1)\cap V$. Let $(u_0 , v_0) \in V$ be the point corresponding to the
complete structure; that is, $\iota^{-1} G^{-1} \{(m_1 ^0, l_1 ^0,m_2 ^0, l_2 ^0) \}$. 

The steps of the proof are as follows: 
\begin{itemize}
\item First, we will show the smoothness of the geometric component $C_1$
as seen in the tetrahedral parameter space where the first cusp 
remains a cusp. This essentially follows from the gradients
of $m_1$ and Theorem \ref{thm:NZ}.
We realize $C_1$ as a graph of a function from one parameter 
to another. 
\item Second, we will show the smoothness of the geometric component 
$C_1^{p_1, q_1}$
in the tetrahedral parameter space where the first cusp has been 
$(p_1, q_1)$-Dehn surgered. The argument is based on continuity 
method for the gradients. We also realize $C_1^{p_1, q_1}$ 
as a graph of a function.
\item We show that the sequence of coefficients of the Taylor series of 
the second functions converges to those of the first function.
\item We change variable to $(m_1, l_1)$.
This change of variables gives us the desired conclusion.
\end{itemize}

(1) The first step is to show the smoothness of $C_1$:

For $(u,v) \in V$, we will denote $G_1 \iota(u,v)$ and $G_2 \iota (u,v)$ by 
$(m_1 , l_1)$ and $(m_2 , l_2)$, considered as functions of $(u,v)$,
respectively. By results of Neumann-Zagier(\cite{NZ}), 
gradient vectors $\nabla m_1$ and $\nabla m_2 $ are linearly independent at 
$(u_0 , v_0)$.($\nabla$ is taken with respect to $(u,v)$). 
Similarly $\nabla l_1$ and $\nabla l_2 $ are linearly independent at $(u_0 , v_0)$.

Moreover, by Theorem \ref{thm:NZ}, $\nabla m_1 (u_0 , v_0)$
and $\nabla l_1 (u_0 , v_0)$ are not linearly dependent over $\bR$
(though they are dependent over $\bC$). So, by the change of
parameter $(u,v) \mapsto (v,u)$ if necessary, we may assume that
 \begin{center}{ \it $\frac{\partial m_1}{\partial v}\left|_{(u_0,v_0)} \right.$ and
$\frac{\partial l_1}{\partial v}\left|_{(u_0,v_0)} \right.$ are linearly
independent over $\bR$. } \end{center} 

Since $\nabla m_1|_{(u_0 , v_0)} \ne 0$, it follows that $\cC_1$ and $C_1$ are 
smooth curves near $\mathbf{z} ^0 $ and $(u_0 , v_0 )$, respectively. 

In addition,
since $\frac{\partial m_1}{\partial v} \ne 0$, $C_1$ is the graph
of a holomorphic map
$$v-v_0 = c_1(u - u_0) +\frac{c_2}{2}(u - u_0)^2 
+ \frac{c_3}{6}(u-u_0)^3 + O((u - u_0)^4)$$ 
at $(u_0 , v_0)$. 
Here $c_1 \ne 0$ since $\frac{\partial m_1}{\partial u}\ne 0$. 

(2) The second step is to do the same for $C_1^{p_1, q_1}$:

For each coprime pair $p_1, q_1$ of
integers for which $|p_1|+|q_1|$ is sufficiently large, 
take $(u_0^{p_1,q_1}, v_0^{p_1 , q_1}) \in \cV$ such that 
\[p_1 \log(\frac{m_1}{m_1 ^0}) + q_1 \log(\frac{l_1}{l_1 ^0}) = \pi \sqrt{-1}\]
so that
$ m_2(u_0^{p_1,q_1}, v_0^{p_1 , q_1}) = m_2^0$ and  
$l_2(u_0^{p_1,q_1}, v_0^{p_1 , q_1}) = l_2 ^0$. 
The point $(u_0^{p_1, q_1}, v_0^{p_1, q_1}$ can be taken to be 
one realized as a hyperbolic manifold obtained 
by $(p_1, q_1)$-Dehn filling the first cusp and the second cusp 
remaining a cusp. Since a sequence of such manifold 
converges to $W$ as $|p_1| + |q_1| \ra \infty$, 
we see that 
$u_0^{p_1,q_1} \rightarrow u_0$ and 
 $v_0^{p_1 , q_1}\rightarrow v_0$ as 
$|p_1| + |q_1| \rightarrow \infty$.

Let $C _1 ^{p_1 , q_1} $ be the set of
points $(u,v)$ in $V$ satisfying 
$p_1 \log(\frac{m_1}{m_1 ^0}) + q_1 \log(\frac{l_1}{l_1 ^0}) = \pi \sqrt{-1}$.

\begin{prop}\label{prop:propA}
When $|p_1| + |q_1|$ is sufficiently large, $C _1^{p_1 , q_1} $
is a smooth curve near
$(u_0^{p_1,q_1}, v_0^{p_1 , q_1})$.
\end{prop}
\begin{proof}
We have that 
$$\frac{\partial}{\partial v}\{p_1\log(\frac{m_1}{m_1 ^0}) 
+ q_1\log(\frac{l_1}{l_1 ^0}) \} 
= \frac{p_1}{m_1}\frac{\partial m_1}{ \partial v} 
+ \frac{q_1}{l_1} \frac{\partial l_1}{ \partial v}.$$ 

Since the values of two functions $\frac{\partial m_1}{ m_1 \partial v}$ and
$\frac{\partial l_1}{ l_1 \partial v}$
are linearly independent over $\bR$ at any point $(u,v) \in V$ near
$(u_0 , v_0)$, $\frac{\partial}{\partial v}\{p_1 \log(\frac{m_1}{m_1 ^0}) + q_1
\log(\frac{l_1}{l_1 ^0}) \}$ is never 0 near $(u_0 , v_0)$. Thus $\nabla\{p_1
\log(\frac{m_1}{m_1 ^0}) + q_1
\log(\frac{l_1}{l_1 ^0})\}$ never vanishes near $(u_0 , v_0)$.
\end{proof}


For $(u,v) \in V$, $m_1 = m_1^0 (= \pm 1)$ if and only
if $l_1 = l_1 ^0$ since parabolic elements of $\slc$ does not
commute with elliptic or loxodromic elements.
Thus $\nabla m_1$ and $\nabla l_1$ are parallel at $(u_0 , v_0)$;
i.e., \[\frac{\partial m_1}{\partial u} \frac{\partial l_1}{\partial v}
- \frac{\partial m_1}{\partial v} \frac{\partial l_1}{\partial u} = 0 
\hbox{ at } (u_0 , v_0).\]

\begin{prop} \label{prop:propB}
The direction of $\nabla \{p_1 \log(\frac{m_1}{m_1 ^0}) + q_1
\log(\frac{l_2}{l_2 ^0})\}(u_0^{p_1,q_1}, v_0^{p_1 , q_1})$ converges to that
of $\nabla m_1(u_0 , v_0)$ as $|p_1| + |q_1| \rightarrow \infty$.
\end{prop}
\begin{proof}
\begin{eqnarray*}
\nabla m_1 = 
({\frac{\partial m_1}{\partial u}}, {\frac{\partial m_1 }{\partial v}})
\end{eqnarray*}
and \begin{eqnarray*}\nabla \left\{p_1 \log(\frac{m_1}{m_1 ^0}) + q_1
\log(\frac{l_2}{l_2 ^0})\right\} 
= \left(\frac{p_1}{m_1} \frac{\partial m_1}{ \partial u}
+ \frac{q_1}{l_1} \frac{\partial l_1}{ \partial u}, 
{\frac{p_1}{m_1} \frac{\partial m_1}{ \partial v} +
\frac{q_1}{l_1}\frac{\partial l_1}{ \partial v}}\right)\end{eqnarray*}
by definition. But
\begin{eqnarray*}
\left. \frac{\frac{p_1}{m_1} \frac{\partial m_1}{ \partial u} +
\frac{q_1}{l_1}\frac{\partial l_1}{ \partial u}}
{\frac{p_1}{m_1} \frac{\partial m_1}{ \partial v}
+ \frac{q_1}{l_1} \frac{\partial l_1}{ \partial v}}
\right|_{(x_0^{p_1,q_1}, y_0^{p_1 , q_1})}  
\rightarrow \left. \frac{\frac{\partial m_1 }{\partial u}}
{\frac{\partial m_1 }{\partial v}}\right|_{(x_0 , y_0)}
\end{eqnarray*}  as $|p_1| + |q_1| \rightarrow
\infty$ by the following  lemma.
\end{proof}

 \begin{lem}\label{lem:pqconv} 
Let $a, b, c, d$ be complex numbers such that $c, d$ are
linearly independent over $\bR$ and $ad - bc=0$. Let $A, B, C, D$ 
be functions on $\bZ \times \bZ$ such that 
\[A(p,q) \rightarrow a, \; B(p,q) \rightarrow b, \; C(p,q) \rightarrow c, \;D(p,q)
 \rightarrow d \hbox{ as } |p| + |q| \rightarrow \infty\] 
Then 
\[\frac{Ap + Bq}{Cp + Dq} \rightarrow \frac{ap+bq}{cp+dq} (= a/c =b/d)
\hbox{ as } |p| + |q| \rightarrow \infty.\]
\end{lem}

(3) The third step is to show the convergence of coefficients.

From Propositions \ref{prop:propA} and \ref{prop:propB}, we see that 
when $|p_1| + |q_1|$ is sufficiently large,
$C_1^{p_1, q_1}$ is the graph of a holomorphic map
\begin{equation*}
\begin{split}
v^{p_1 , q_1} - v_0^{p_1,q_1}= c_1^{p_1,q_1} (u^{p_1,q_1} - u_0^{p_1,q_1}) + 
\frac{c_2^{p_1 , q_1}}{2} (u^{p_1,q_1} - u_0^{p_1,q_1})^2 \\
+\frac{ c_3^{p_1 , q_1} }{6} (u^{p_1,q_1} - u_0^{p_1,q_1})^3 
+ O((u^{p_1,q_1} - u_0^{p_1,q_1})^4)
\end{split}
\end{equation*}
at $(u_0 ^{p_1,q_1}, v_0 ^{p_1,q_1}) $, where $c_1 ^{p_1,q_1} \ne 0$.

Let $x$ denote the variable on $C_1$ by setting it equal to $u$ for 
a point $(u, v)$ on $C_1$. Let $y$ denote the function on 
$C_1$ set to be equal to $v$. We regard $y$ as a function of $x$. 

\begin{prop}\label{prop:convci}
$c_i ^{p_1 , q_1} \rightarrow c_i $ 
as $|p_1| + |q_1| \rightarrow \infty$ for $i=1 ,2, 3$.
\end{prop}
\begin{proof}
Any $(u,v) \in C_1$ satisfies
\begin{equation*}
\begin{split}
p_1 \left\{\frac{1}{m_1 ^0} (m_1 - m_1 ^0 ) - \frac{1}{2} (m_1-m_1 ^0)^2  
+ \frac{1}{3 m_1 ^0}(m_1 -m_1 ^0)^3 \right\} \\+ q_1
\left\{ \frac{1}{l_1 ^0} (l_1 -l_1 ^0) - \frac{1}{2} (l_1 -l_1^0)^2 
+ \frac{1}{3 l_1 ^0}(l_1 -l_1 ^0)^3  \right\} = 0
\end{split}
\end{equation*}
up to order three for any pair $p_1$,$q_1$ of integers. 
The left hand side of the
equation is the Taylor expansion of the function 
$p_1\log(\frac{m_1}{m_1^0}) + q_1 \log(\frac{l_1}{l_1^0})$ 
at $(m_1^0 , l_1 ^0)$ up to order 3 terms. 

Differentiating above equation
with respect to $x$ gives us the relation up to order three
\begin{equation} \label{eq:d1}
\begin{split}
p_1 \left\{ \frac{1}{m_1^0} \frac{dm_1}{dx} - (m_1 - m_1^0) \frac{dm_1}{dx}
+ \frac{1}{m_1^0} (m_1 - m_1^0)^2 \frac{dm_1}{dx} \right\}
\\ + q_1 \left\{ \frac{1}{l_1^0} \frac{dl_1}{dx} - (l_1 - l_1^0) \frac{dl_1}{dx}
+ \frac{1}{l_1^0} (l_1 - l_1^0) ^2 \frac{dl_1}{dx} \right\} = 0.
\end{split}
\end{equation}

 Applying the chain rule 
$$\frac{df}{dx} = \frac{\partial f}{\partial u} +
\frac{dy}{dx}\;\frac{\partial f}{\partial v}$$ to $f=m_1$ 
and $f=l_1$ in (\ref{eq:d1}), 
we obtain
$$c_1 = 
\left. \frac{dy}{dx}\right|_{x=u_0}  
= \left. -\frac{ \frac{p_1}{m_1^0} \frac{\partial m_1}{\partial u} +  
\frac{q_1}{l_1^0}\frac{\partial l_1}{\partial u}}
{\frac{p_1}{m_1^0}\frac{\partial m_1}{\partial v} +  \frac{q_1}{l_1^0}
\frac{\partial l_1}{\partial v} }\right|_{(u,v)=(x_0 , y_0)}.
$$
On the other hand, any $(u,v) \in C_1^{p_1,q_1}$ satisfies
$$p_1 \log(\frac{m_1}{m_1^0}) + q_1 \log(\frac{l_1}{l_1^0}) = \pi \sqrt{-1} . $$

We let $x^{p_1, q_1}$ denote the variable on $C_1^{p_1, q_1}$ 
obtained by setting it equal to $u$ and 
let $y^{p_1, q_1}$ denote the function on the same by setting it equal to $v$. 
We consider $y^{p_1, q_1}$ the function of $x^{p_1, q_1}$. 

Differentiating with respect to $x^{p_1 , q_1 }$ gives us the relation
\begin{equation}\label{eq:d1pq}
\frac{p_1}{m_1} \frac{dm_1}{dx^{p_1 , q_1 }} + \frac{q_1}{l_1} \frac{dl_1}{dx^{p_1 , q_1 }} =0
\end{equation}
Then using the chain rule in (\ref{eq:d1pq}), we obtain
$$c_1 ^{p_1,q_1} 
= \left. \frac{dy^{p_1,q_1}}{dx ^{p_1 , q_1 }}\right|_{x^{p_1 , q_1 }= x_0^{p_1 , q_1 }}
= \left. -\frac{\frac{p_1}{m_1} \frac{\partial m_1}{ \partial u} +
\frac{q_1}{l_1}\frac{\partial l_1}{ \partial u}}{\frac{p_1}{m_1} 
\frac{\partial m_1}{ \partial v} + 
\frac{q_1}{l_1} \frac{\partial l_1}{ \partial v}}
\right|_{(u,v) = (u_0^{p_1,q_1}, v_0^{p_1 , q_1})}.
$$
From Lemma \ref{lem:pqconv}, we see that 
$c_1^{p_1,q_1}\rightarrow  c_1$ as $|p_1| + |q_1| \rightarrow \infty$.

Differentiating (\ref{eq:d1}) with respect to $x$ and using
$$\frac{d^2 f}{dx^2} = \frac{\partial ^2 f}{\partial u^2}
+ 2 \frac{dy}{dx}\frac{\partial ^2 f}{\partial u \partial v}
+ ( \frac{dy}{dx}) ^2 \frac{\partial ^2 f}{\partial v^2}
+ \frac{d^2 y}{dx^2}\frac{\partial f}{\partial v},$$
we obtain
\begin{equation*}
\begin{split}
c_2 &=\left. \frac{ d^2 y}{dx^2}\right|_{x=u_0} \\
=& -[ \frac{p_1}{m_1 ^0} \{\frac{\partial ^2 m_1}{\partial u^2}
+ 2 c_1 \frac{\partial ^2 m_1}{\partial u \partial v}
+ c_1 ^2 \frac{\partial ^2 m_1}{\partial v^2}
- \frac{1}{m_1 ^0} ( \frac{\partial m_1}{\partial u} + c_1 \frac{\partial m_1}{\partial v})^2
\}\\
+& \frac{q_1}{l_1 ^0} \{\frac{\partial ^2 l_1}{\partial u^2}
+ 2 c_1 \frac{\partial ^2 l_1}{\partial u \partial v}
+ c_1 ^2 \frac{\partial ^2 l_1}{\partial v^2}
- \frac{1}{l_1 ^0} ( \frac{\partial l_1}{\partial u} + c_1\frac{\partial l_1}{\partial v} )^2 \}
]\\
&/\left. \left[\frac{p_1}{m_1 ^0}\frac{\partial m_1}{\partial v} +  
\frac{q_1}{l_1 ^0}\frac{\partial l_1}{\partial v}\right]
\right|_{(u,v)=(u_0 , v_0 )}
\end{split}
\end{equation*}

But differentiating (\ref{eq:d1pq}) with respect to $x^{p_1 , q_1 }$, we obtain
\begin{equation*}
\begin{split}
c_2^{p_1, q_1} &= 
\frac{ d^2 y^{p_1 , q_1}}{d(x ^{p_1 , q_1})^2}\left|_{x ^{p_1 , q_1 }=x_0^{p_1 , q_1}} \right. \\
=& - [ \frac{p_1}{m_1}
\{\frac{\partial ^2 m_1}{\partial u^2}
+ 2 c_1 ^{p_1, q_1} \frac{\partial ^2 m_1}{\partial u \partial v}
+ (c_1 ^{p_1, q_1}) ^2 \frac{\partial ^2 m_1}{\partial v^2}
- \frac{1}{m_1} ( \frac{\partial m_1}{\partial u} 
+ c_1 ^{p_1, q_1} \frac{\partial m_1}{\partial v} )^2 \}
\\+ & \frac{q_1}{l_1} \{\frac{\partial ^2 l_1}{\partial u^2}
+ 2 c_1 ^{p_1, q_1} \frac{\partial ^2 l_1}{\partial u \partial v}
+ (c_1 ^{p_1, q_1}) ^2 \frac{\partial ^2 l_1}{\partial v^2}
- \frac{1}{l_1} ( \frac{\partial l_1}{\partial u} 
+ c_1 ^{p_1, q_1} \frac{\partial l_1}{\partial v} )^2 \}  ]\\
& /
[\frac{p_1}{m_1} \frac{\partial m_1}{\partial v} 
+ \frac{q_1}{l_1} \frac{\partial l_1}{\partial v}]
\left|_{(u,v)=(x_0 ^{p_1 , q_1} , y_0 ^{p_1 , q_1})}\right.
\end{split}
\end{equation*}
Using Lemma \ref{lem:pqconv} again, we see that
$c_2^{p_1,q_1}\rightarrow c_2$ as $|p_1| + |q_1| \rightarrow \infty$. 
Continuing in this way, we can also show that  
$c_3^{p_1,q_1}\rightarrow  c_3$ as $|p_1| + |q_1| \rightarrow \infty$.
\end{proof}

\begin{figure}
\begin{center}
\epsfig{file=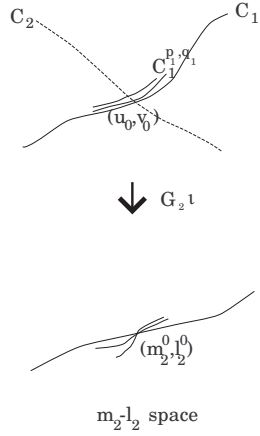 , height = 6cm}
\end{center}
\caption{Convergence of smooth curves} \label{fig:convcurves}
\end{figure}


(4) The fourth step is to consider the holomorphic map 
$G_2 \iota = (m_2 , l_2 ): V \rightarrow \bC^2$ defined near $(x_0 , y_0)$.

\begin{prop} The images $\cC_1$ and $\cC_1^{p_1, q_1}$ of 
$C_1$ and $C_1^{p_1 , q_1}$, when $|p_1| + |q_1|$ is sufficiently large,
under $G_2 \iota$ respectively are
smooth curves through $(m_2^0 , l_2^0)$. $\cC_1$ is the graph of a holomorphic map
\begin{equation} \label{eq:linminf}
\begin{split}
l_2 = l_2 ^0 + a_1 (m_2 - m_2 ^0) + \frac{a_2}{2} (m_2 - m_2 ^0)^2+ \frac{a_3}{6} (m_2 - m_2 ^0 )^3 \\
+ \mbox{higher order terms}.
\end{split}
\end{equation}
near $(m_2 ^0 , l_2 ^0)$ with $a_1 \ne 0$ and 
$\cC_1^{p_1 , q_1 }$ is the graph of a holomorphic map
\begin{equation} \label{eq:linmpq}
\begin{split}
l_2 = l_2^0 + a_1^{p_1 , q_1} (m_2 - m_2^0) + \frac{a_2^{p_1 , q_1}}{2} (m_2 - m_2^0)^2
+ \frac{a_3^{p_1 , q_1}}{6} (m_2 - m_2^0 )^3 \\
+ \mbox{higher order terms}.
\end{split}
\end{equation}
\end{prop}

\begin{proof}
$C_1$ is the graph of a holomorphic map with a Taylor expansion
$$ v - v_0 = c_1 (u - u_0) + \frac{c_2}{2} (u - u_0)^2 + 
\frac{c_3}{6}(u-u_0)^3 + O((u - u_0)^4) $$ near $(u_0 , v_0)$. 
By our definition of $C_1$, we
have $m_1 = m_1^0$ for points on $C_1$. Thus
$\frac{dm_1}{dx}\left|_{x= u_0} \right. = [\frac{\partial m_1}{\partial u} + 
c_1\frac{\partial m_1}{\partial v}]\left|_{(u_0 , v_0) } = 0 \right.$. 

Since $\nabla m_1 (u_0 , v_0)$ and $\nabla m_2 (u_0 , v_0)$ are linearly
independent, we must have 
\[\left. \frac{dm_2}{dx}\right|_{x= u_0} =
\left. \left[\frac{\partial m_2}{\partial u} + 
c_1 \frac{\partial m_2}{\partial v}\right]\right|_{(u,v)=(u_0 , v_0) } \ne 0.\] 
Thus, the map $G_2 \iota$ restricted to $C_1$ is nonsingular at $(u_0 , v_0)$;
hence, the image of $C_1$ under $G_2 \iota$ is a smooth curve
through $(m_2^0 , l_2^0)$.  

Recall that $C_1 ^{p_1 , q_1}$ is
the graph of a holomorphic map with an expansion
\begin{equation*}
\begin{split}
v - v_0^{p_1,q_1}= c_1 ^{p_1,q_1} (u - u_0^{p_1,q_1}) 
+ \frac{c_2 ^{p_1 , q_1}}{2} (u - u_0 ^{p_1,q_1})^2 \\
+ \frac{c_3 ^{p_1 , q_1}}{6} (u - u_0 ^{p_1,q_1})^3 +
O((u - u_0 ^{p_1, q_1})^4)
\end{split}
\end{equation*}
at $(u_0^{p_1,q_1}, v_0 ^{p_1,q_1}) $. 

Since $c_1^{p_1 , q_1} \rightarrow c_1$ 
and $(u_0 ^{p_1,q_1}, v_0 ^{p_1,q_1}) \rightarrow (u_0 , v_0)$ 
as $|p_1| + |q_1| \rightarrow \infty$, we have
\[\left. \frac{dm_2}{dx}\right|_{x = u_0 ^{p_1,q_1} }
=\left.\left[\frac{\partial m_2}{\partial u} + c_1 ^{p_1 , q_1}
\frac{\partial m_2}{\partial v}\right]\right|_{(u,v)
=(u_0^{p_1 , q_1}, v_0^{p_1 , q_1}) } \ne 0\] 
when $|p_1| + |q_1|$ is sufficiently large.

We showed above 
\[\left.\frac{dm_2}{dx}\right|_{x= x_0} \ne 0 \hbox{ and }  
\left. \frac{dm_2}{dx^{p_1 , q_1}}\right|_{x^{p_1 , q_1} = x_0^{p_1,q_1} } \ne 0\] 
when $|p_1| + |q_1|$ is sufficiently large. In
the same way we show 
\[\left.\frac{dl_2}{dx}\right|_{x= u_0} \ne 0 \hbox{ and } 
\left. \frac{dl_2}{dx}\right|_{x = u_0 ^{p_1,q_1} } \ne 0\] 
when $|p_1| + |q_1|$ is sufficiently large. 
Thus, the rest of the proposition follow.
\end{proof}
 
\begin{proof}[Proof of Theorem \ref{thm:convcoeff}] 
We differentiate (\ref{eq:linminf}) with respect to $x$ and evaluate
at $x = u_0$ to obtain
$$ \left. \frac{dl_1}{dx}\right|_{x = u_0}  
= a_1 \left. \frac{dm_1}{dx} \right|_{x = u_0 }  .$$
So 
\[a_1 = \left. \frac{\frac{\partial l_1}{\partial u} + c_1 \frac{\partial l_1}{\partial v} }
{\frac{\partial m_1}{\partial u} + c_1 \frac{\partial m_1}{\partial v}}\right|_{(u,v)=(u_0 , v_0)} .
\] 
Similarly, from (\ref{eq:linmpq}), we obtain
$$a_1 ^{p_1 , q_1 }  
= \left. \frac{\frac{\partial l_1}{\partial u} + c_1 ^{p_1 , q_1 }\frac{\partial l_1}{\partial v} }
{\frac{\partial m_1}{\partial u} + c_1^{p_1 , q_1 }\frac{\partial m_1}{\partial v}}
\right|_{(u,v)=(x_0 ^{p_1 , q_1 }, y_0^{p_1 , q_1 })} .$$ 
Since $c_1^{p_1 , q_1} \rightarrow c_1 $ as $|p_1| + |q_1| \rightarrow \infty $ by
Proposition \ref{prop:convci}, 
we have  $a_1^{p_1 , q_1} \rightarrow a_1 $ as $|p_1| + |q_1| \rightarrow\infty $. 

Proceeding in this way using successive
differentiations we also obtain $a_i^{p_1 , q_1} \rightarrow a_i $ as 
$|p_1| + |q_1| \rightarrow \infty $ for $i = 2,3$.
\end{proof}

We close this section with the following lemma which is implied
by the fact that
the curves $G_2 \iota (C_1 ^{p_1 , q_1})$ are invariant under the
involution $(m_2 , l_2) \mapsto (1 / m_2 , 1 / l_2 )$ restricted
to a small neighborhood of $(m_2 ^0 , l_2 ^0)$ in $\bC^2$.

\begin{lem} \label{lem:involution}
When $|p_1| + |q_1|$ is large, 
$a_2^{p_1 , q_1} = -m_2^0 a_1^{p_1 , q_1} + l_2^0 (a_1^{p_1 , q_1})^2$.
\end{lem}
\begin{proof} 
By results of Neumann and Zagier(\cite{NZ}),
there is an involution $\mathcal{I}: V \rightarrow V$ such that
$G_1 \iota \mathcal{I} (u,v) = (m_1 , l_1)$ and
$G_2 \iota \mathcal{I} (u,v) = (1 / m_2 , 1 / l_2)$ when $G_1
\iota (u,v) = (m_1 , l_1)$ and $G_2
\iota (u,v) = (m_2 , l_2)$. Since $C_1 ^{p_1 , q_1}$ is invariant
under $\mathcal{I}$, $G_2 \iota (C_1 ^{p_1 , q_1})$ is invariant
under the involution $(m_2 , l_2) \mapsto (1 / m_2 , 1 / l_2)$.
\end{proof}

\end{document}